\documentclass[10pt]{amsart}%{scrartcl,scrreprt}
\usepackage[T1]{fontenc}
\usepackage{lmodern}
\usepackage{geometry}
\usepackage[latin1] {inputenc}
\usepackage{amsmath}
\usepackage{amsfonts, amssymb, textcomp}
\usepackage[colorlinks=flase, linkcolor=red,urlcolor=green, citecolor=blue]{hyperref}
\usepackage{subeqnarray}

\usepackage{latexsym}
\usepackage{fancyhdr}
\usepackage{longtable}
\usepackage{amsmath, amssymb}
\usepackage{graphicx}

\numberwithin{equation}{section}

\theoremstyle{plain}
\newtheorem{theorem}{Theorem}[section]
\newtheorem{lemma}[theorem]{Lemma}

\newtheorem{corollary}[theorem]{Corollary}
\theoremstyle{definition}
\newtheorem{definition}[theorem]{Definition}
\newtheorem{example}[theorem]{Example}
\newtheorem{remark}[theorem]{Remark}

\newcommand{\RR}{\mathbb{R}}
\newcommand{\CC}{\mathbb{C}}
\newcommand{\NN}{\mathbb{N}}

\newcommand{\DD}{\mathbb{D}}

\newcommand{\id}{\operatorname{id}}

\let\on=\operatorname
\newenvironment{demo}[1]{\par\smallskip\noindent{\bf #1.}}{\par\smallskip}

\makeatletter
\@namedef{subjclassname@2020}{%
  \textup{2020} Mathematics Subject Classification}
\makeatother

\title[Operators on weighted analytic spaces]
{The differentiation and integration operator on weighted spaces of analytic functions in the weight sequence setting}

\author[G.~Schindl]{Gerhard Schindl}

\address{G.~Schindl: Fakult\"at f\"ur Mathematik, Universit\"at Wien, Oskar-Morgenstern-Platz~1, A-1090 Wien, Austria.}
\email{gerhard.schindl@univie.ac.at}

\begin{document}

\begin{abstract}
We investigate the differentiation and integration operator defined on resp. between weighted spaces of analytic functions in the special situation when the weight is given in terms of a weight sequence. More precisely, in this case the weight is defined via the associated weight function which is a crucial object in the ultradifferentiable setting and the information can be then expressed purely in terms of the underlying weight sequence. Due to the additional structure this setting is exploited in order to construct (counter-)examples, to compare conditions on weights in the different weighted frameworks, and to study and to understand in a more precise qualitative way required and admissible growth restrictions on the weights.
\end{abstract}

\thanks{This research was funded in whole or in part by the Austrian Science Fund (FWF) 10.55776/PAT9445424}
\date{\today}

\maketitle

\section{Introduction}
The aim of this article is the continuation of the author's study from \cite{solidassociatedweight}, \cite{weightedentireinclusion1} and \cite{weightedentireinclusion2} dealing with certain weighted spaces of analytic functions. Let $v:[0,+\infty)\rightarrow(0,+\infty)$ be the \emph{strictly positive and radial weight (function)} and assume that $v$ is
\begin{itemize}
\item[$(a)$] continuous,
	
\item[$(b)$] non-increasing and

\item[$(c)$] rapidly decreasing, i.e. $\lim_{t\rightarrow+\infty}t^kv(t)=0$ for all $k\ge 0$ (resp. equivalently for all $k\in\NN$).
\end{itemize}
This is the same notion to be a weight (function) as it has been considered in \cite{Bonet2022survey}, see \cite[Sect. 2, p. 3]{Bonet2022survey}, and also in \cite{solidassociatedweight}, \cite{weightedentireinclusion1} and \cite{weightedentireinclusion2}; and $v$ is called \emph{normalized} if $v(t)=1$ for all $t\in[0,1]$. Note that in the literature there exist different notions to be a weight function in the weighted analytic setting; occasionally one deals with more general weights (i.e. weaker assumptions von $v$). Moreover, some authors use the notation $v^{-1}$ instead of $v$, see e.g. \cite{AbaninThien17}, and, indeed, this comment is crucial when comparing (characterizing) conditions for weight functions since in general the inversion $v\mapsto v^{-1}$ does not preserve growth requirements.\vspace{6pt}

The corresponding \emph{weighted (Banach) space of entire functions} is defined as follows:
\begin{equation}\label{weightedclassdef}
H^{\infty}_v(\CC):=\{f\in H(\CC): \|f\|_v:=\sup_{z\in\CC}|f(z)|v(|z|)<+\infty\}.
\end{equation}
$H(\CC)$ denotes the class of entire functions and analogously in the literature spaces with $o$-growth restriction are considered:
$$H^0_v(\CC):=\{f\in H(\CC): \lim_{|z|\rightarrow+\infty}|f(z)|v(|z|)=0\};$$
see e.g. \cite[Sect. 2]{Bonet2022survey}.

Recall that normalization can be assumed w.l.o.g. since for each $v$ one can find a normalized $v^n$ such that $v(z)=v^n(z)$ for all $z\in\CC$ with $|z|$ large enough and hence $H^{\infty}_v(\CC)=H^{\infty}_{v^n}(\CC)$ resp. $H^0_v(\CC)=H^0_{v^n}(\CC)$ as Banach spaces; see Remark \ref{normrem} and Lemmas \ref{normconvexrem}, \ref{auxlemma}. Frequently one considers also the \emph{unit disk case,} write $H^{\infty}_{v_{\mathbb{D}}}(\mathbb{D})$ resp. $H^0_{v_{\mathbb{D}}}(\mathbb{D})$, and here one assumes that $v_{\mathbb{D}}: [0,1)\rightarrow(0,+\infty)$ satisfies $(a)$ and $(b)$, whereas $(c)$ is replaced by $\lim_{t\rightarrow 1^{-}}t^kv(t)=0$ for all $k\ge 0$; see again \cite[Sect. 2, p. 3]{Bonet2022survey}. Indeed, this last property is equivalent to $\lim_{t\rightarrow 1^{-}}v(t)=0$. Accordingly, we say that $v_{\mathbb{D}}$ is normalized if $v_{\mathbb{D}}(0)=1$ and which can analogously be assumed w.l.o.g. (see again Remark \ref{normrem}).

For our considerations in this present work, in the weighted entire case, we focus on $v=v_{\mathbf{M}}$ with
\begin{equation}\label{weights}
v_{\mathbf{M}}(t):=\exp(-\omega_{\mathbf{M}}(t)),\;\;\;t\ge 0.
\end{equation}
$\mathbf{M}\in\RR_{>0}^{\NN}$ denotes a sequence of positive real numbers satisfying basic growth conditions (see Definition \ref{defweightsequ}) and $\omega_{\mathbf{M}}$ is the corresponding \emph{associated weight function;} we refer to Section \ref{assofunctionsection} for the precise definition, more details and citations. The basic growth properties of $\omega_{\mathbf{M}}$ ensure that $v_{\mathbf{M}}$ is a weight function in the above sense and we have investigated this case in \cite{solidassociatedweight} concerning the study of solid hulls and cores, whereas in \cite{weightedentireinclusion1} and \cite{weightedentireinclusion2} the inclusion relations have been treated in detail. The main motivation for this approach is to exploit the additional information and to characterize desired properties of the weighted spaces $H^{\infty}_{v_{\mathbf{M}}}(\CC)$ and $H^{0}_{v_{\mathbf{M}}}(\CC)$ in terms of $\mathbf{M}$ which is reasonable since $v_{\mathbf{M}}$ is purely expressed in terms of the weight sequence $\mathbf{M}$. Thus this setting is also tailor-made to study and construct explicit (counter-)examples. Moreover, the second central idea is to connect and to compare different areas of research dealing with weighted spaces. For this recall that $\omega_{\mathbf{M}}$ plays a prominent and crucial role in the theory of \emph{ultradifferentiable and ultraholomorphic function classes,} see e.g. \cite{Komatsu73}, \cite{BonetMeiseMelikhov07}, \cite{optimalflat23}, and one is interested in comparing the appearing conditions for $\mathbf{M}$ and to study the meaning of known properties for $\mathbf{M}$ in the weighted entire setting as well.

Moreover, recall that for any abstractly given weight function $v$ one can associate in a straight-forward way a sequence $\mathbf{M}^v$ satisfying standard regularity conditions, see Section \ref{assosequencesection} and the references there. Hence there exists a natural connection to the weight sequence setting and one asks how properties can be transferred between spaces defined via $v$ and in terms of $\mathbf{M}^v$.\vspace{6pt}

Similarly, one considers the \emph{weighted unit disk case,} see \cite[Sect. 7]{solidassociatedweight} (with $c=1$). We set
\begin{equation}\label{weightfctunitdisk}
v_{\mathbf{M},\mathbb{D}}(t)=\exp(-\omega_{\mathbf{M},\mathbb{D}}(t)):=\exp\left(-\omega_{\mathbf{M}}\left(\frac{1}{1-t}\right)\right),\;\;\;t\in[0,1),
\end{equation}
and write $H^{\infty}_{v_{\mathbf{M},\mathbb{D}}}(\DD)$ resp. $H^{0}_{v_{\mathbf{M},\mathbb{D}}}(\DD)$ for the corresponding weighted space.\vspace{6pt}

The survey article \cite{Bonet2022survey} provides a recent, detailed and exhaustive summary of statements concerning weighted spaces of analytic functions defined on $\CC$ and $\DD$ and focuses also on various operators acting on resp. between these weighted spaces. Motivated by the results gathered in \cite[Sect. 8]{Bonet2022survey} and the approach in \cite{AbaninThien17} we investigate in this work the \emph{differentiation operator $D$} and the \emph{integration operator $J$} in the weight sequence setting; i.e. we focus on properties of $D$ and $J$ when acting on and between $H^{\infty}_{v_{\mathbf{M}}}(\CC)$, $H^{\infty}_{v_{\mathbf{M},\mathbb{D}}}(\DD)$, $H^{0}_{v_{\mathbf{M}}}(\CC)$, $H^{0}_{v_{\mathbf{M},\mathbb{D}}}(\DD)$. Note that the weight $v$ in \cite{AbaninThien17} corresponds to $v^{-1}$ in \cite{Bonet2022survey}.

We provide information purely expressed in terms of the defining sequence $\mathbf{M}$ and relate the conditions to known properties appearing in the ultradifferentiable setting. In \cite[Sect. 4]{weightedentireinclusion1}, see also the correction in \cite[Appendix A]{weightedentireinclusion2}, we have already investigated the \emph{bilinear multiplication operator} $\mathfrak{m}: (f,g)\mapsto f\cdot g$. It turns out that in general $H^{\infty}_v(\CC)$ is not stable under this operator, see \cite[Prop. 4.2 \& Cor. 4.14]{weightedentireinclusion1} and the erratum in \cite[Appendix A]{weightedentireinclusion2}. The situation is different for $H^{\infty}_{v_{\mathbb{D}}}(\mathbb{D})$ (see \cite[Prop. 2.1]{BonetDomanskiLindstroem99}) and in the weighted entire setting when considering \emph{weight systems} subject to a parameter $c>0$ and which allows for more flexibility; see the main results \cite[Prop. 4.1, Thm. 4.8 \& 4.13]{weightedentireinclusion1} and also Section \ref{parametersection} for definitions. The plan is to continue the recent study for the \emph{weighted composition operator} in a separate work.\vspace{6pt}

This paper is structured as follows: In Section \ref{weightsection} we gather and recall the relevant notation and conditions for weight functions and weight sequences. By collecting some preparatory material we are able to show the first main Theorem \ref{Thm35disk} which provides the characterization of the inclusion relation $H^{\infty}_{v_{\mathbf{M},\mathbb{D}}}(\DD)\subseteq H^{\infty}_{v_{\mathbf{N},\mathbb{D}}}(\DD)$ resp. $H^{0}_{v_{\mathbf{M},\mathbb{D}}}(\DD)\subseteq H^{0}_{v_{\mathbf{N},\mathbb{D}}}(\DD)$ in terms of condition
\begin{equation}\label{equ14}
\exists\;A\ge 1\;\forall\;p\in\NN:\;\;\;N_p\le AM_p,
\end{equation}
and hence transfers crucial results from \cite{weightedentireinclusion1} and \cite{weightedentireinclusion2} to the \emph{unit disk case.} Then, in Section \ref{diffintsection} we investigate properties for the operators $D$ and $J$ and split the information. Section \ref{diffintsectionOgrowth} is dedicated to the study of continuity, nuclearity and surjectivity of $D$ and $J$ in the weighted entire setting with $O$-growth condition; see Theorems \ref{Diffopthm}, \ref{Intopthm}, \ref{Intopthmcor}, \ref{Intopthmcoru}  and Corollaries \ref{Diffopthmcor}, \ref{DiffIntopthmcor}, \ref{finalDiffIntcor}. In Section \ref{diffintsectionogrowth} we focus on weighted spaces of entire functions with $o$-growth condition and investigate dynamics of $D$ and $J$ on such spaces; see Theorems \ref{Thm33}, \ref{hypercyclthm34}, \ref{hypercyclthm34fctversion} and \ref{hypercyclcor36}. The final Section \ref{unitdiscsection} is dedicated to the study of continuity and surjectivity of $D$ and $J$ in the unit disk case; see Theorem \ref{Thm30&31}. In each section (counter-)examples of weight sequences are constructed in order to illustrate the required growth conditions and to emphasize the weight sequence setting. Summarizing, it turns out that the sequence $\mathbf{G}^1:=(p!)_{p\in\NN}$ resp. the relative growth between the defining sequence $\mathbf{M}$ and $\mathbf{G}^1$ is crucial in many considerations.\vspace{6pt}

\textbf{Acknowledgements.} The author thanks Prof. em. Jos\'{e} Bonet Solves from the Universitat Polit\`{e}cnica de Val\`{e}ncia for helpful and clarifying explanations during the preparation of this work, especially for commenting on the crucial growth conditions of the weights in \cite{AbaninThien17}.

\section{Weight sequences and weight functions}\label{weightsection}
\subsection{Basic notation}
We write $\NN=\{0,1,2,\dots\}$ and $\NN_{>0}=\NN\backslash\{0\}=\{1,2,\dots\}$.

\subsection{Weight sequences}\label{weightsequencesection}
For given $\mathbf{M}=(M_p)_p\in\RR_{>0}^{\NN}$ set $\mu_p:=\frac{M_p}{M_{p-1}}$, $p\in\NN_{>0}$, $\mu_0:=1$. This notation allows the convenient representation $M_p=M_0\prod_{i=1}^p\mu_i$. $\mathbf{M}$ is called \emph{normalized} if $1=M_0\le M_1\Leftrightarrow 1\le\mu_1$ and \emph{log-convex} if
$$\forall\;p\in\NN_{>0}:\;\;\;M_p^2\le M_{p+1}M_{p-1},$$
which is equivalent to the fact that $p\mapsto\mu_p$, $p\in\NN_{>0}$, is non-decreasing. Log-convexity is precisely (M.1) in the classical work \cite{Komatsu73}. Set $m_p:=\frac{M_p}{p!}$, and when even $\mathbf{m}=(m_p)_p$ is log-convex, which is equivalent to $p\mapsto\frac{\mu_p}{p}$ being non-decreasing, we call $\mathbf{M}$ \emph{strongly log-convex.} We consider the following basic set of sequences; see \cite[Def. 2.4]{weightedentireinclusion1}:

\begin{definition}\label{defweightsequ}
$\mathbf{M}=(M_p)_p\in\RR_{>0}^{\NN}$ is called a \emph{weight sequence} if $M_0=1$ and
$$\lim_{p\rightarrow+\infty}(M_p)^{1/p}=+\infty.$$
\end{definition}
Note that for any log-convex weight sequence the mapping $p\mapsto(M_p)^{1/p}$ is non-decreasing; recall e.g. \cite[Lemma 2.0.4]{diploma}.

Occasionally we consider the following (sub)set of weight sequences:
$$\hypertarget{LCset}{\mathcal{LC}}:=\{\mathbf{M}\in\RR_{>0}^{\NN}:\;\mathbf{M}\;\text{is normalized, log-convex},\;\lim_{p\rightarrow+\infty}(M_p)^{1/p}=+\infty\}.$$
We see that $\mathbf{M}\in\hyperlink{LCset}{\mathcal{LC}}$ if and only if $1=\mu_0\le\mu_1\le\dots$ and $\lim_{p\rightarrow+\infty}\mu_p=+\infty$ (see e.g. \cite[p. 104]{compositionpaper}) and there is a one-to-one correspondence between $\mathbf{M}$ and $\mu=(\mu_p)_p$ by taking $M_p:=\prod_{i=0}^p\mu_i$.

$\mathbf{M}$ satisfies \emph{moderate growth,} abbreviated in the following by \hypertarget{mg}{$(\on{mg})$}, if
$$\exists\;C\ge 1\;\forall\;p,q\in\NN:\;M_{p+q}\le C^{p+q+1} M_p M_q.$$
In view of \cite{Komatsu73} this condition is also known under \emph{(M.2) or stability under ultradifferential operators.} If $M_0=1$, then it suffices to consider $C^{p+q}$ in the previous estimate. It is known that for any log-convex weight sequence property \hyperlink{mg}{$(\on{mg})$} is equivalent to the following growth condition:
\begin{equation}\label{quotientrootcompprop}
\exists\;A\ge 1\;\forall\;p\in\NN_{>0}:\;\;\;\mu_p\le A(M_p)^{1/p}.
\end{equation}
For this we refer to \cite[Lemma 2.2 \& Rem. 2.3]{whitneyextensionweightmatrix}, to the detailed discussion in \cite{modgrowthstrange} and the references therein.\vspace{6pt}

Let $\mathbf{M},\mathbf{N}\in\RR_{>0}^{\NN}$ be given, then write
\begin{itemize}
\item[$(*)$] $\mathbf{M}\le\mathbf{N}$ if $M_p\le N_p$ for all $p\in\NN$,
	
\item[$(*)$] $\mathbf{M}\hypertarget{preceq}{\preceq}\mathbf{N}$ if
	$$\sup_{p\in\NN_{>0}}\left(\frac{M_p}{N_p}\right)^{1/p}<+\infty,$$
	
\item[$(*)$] $\mathbf{M}$ and $\mathbf{N}$ are called \emph{equivalent,} formally denoted by $\mathbf{M}\hypertarget{approx}{\approx}\mathbf{N}$, if $\mathbf{M}\hyperlink{preceq}{\preceq}\mathbf{N}$ and $\mathbf{N}\hyperlink{preceq}{\preceq}\mathbf{M}$.

\item[$(*)$] $\mathbf{M}\hypertarget{mtriangle}{\vartriangleleft}\mathbf{N}$ if $\lim_{p\rightarrow+\infty}\left(\frac{M_p}{N_p}\right)^{1/p}=0$; i.e. if
\begin{equation*}\label{triangleestim}
		\forall\;h>0\;\exists\;C_h\ge 1\;\forall\;p\in\NN:\;\;\;M_p\le C_hh^pN_p.
\end{equation*}
	
\item[$(*)$] Obviously, $\mathbf{M}\hyperlink{mtriangle}{\vartriangleleft}\mathbf{N}$ implies $\mathbf{M}\hyperlink{preceq}{\preceq}\mathbf{N}$, but \hyperlink{mtriangle}{$\vartriangleleft$} is never reflexive.
\end{itemize}
A prominent example are the \emph{Gevrey sequences} $\mathbf{G}^s:=(p!^s)_{p\in\NN}$, $s>0$. Each $\mathbf{G}^s$ satisfies \hyperlink{mg}{$(\on{mg})$}, $\mathbf{G}^s\in\hyperlink{LCset}{\mathcal{LC}}$ for any $s>0$ and $\mathbf{G}^s\hyperlink{mtriangle}{\vartriangleleft}\mathbf{G}^t$ for all $0<s<t$. Especially, the sequence $\mathbf{G}^1$ plays a key-role in our considerations in this paper.

The above growth relations are known and crucial in the ultradifferentiable setting; we refer to \cite[Prop. 2.12 \& 4.6]{compositionpaper} and note that \hyperlink{mg}{$(\on{mg})$} is obviously preserved under equivalence. However, in the weighted entire setting an even stronger relation is becoming relevant; see \cite[$(1.4)$, Thm. 3.5]{weightedentireinclusion1}, \cite[Thm. 3.2]{weightedentireinclusion2}: For log-convex weight sequences $\mathbf{M}$, $\mathbf{N}$ it holds that \eqref{equ14} is valid if and only if $H^{\infty}_{v_{\mathbf{M}}}(\CC)\subseteq H^{\infty}_{v_{\mathbf{N}}}(\CC)$ (with continuous inclusion) and if and only if $v_{\mathbf{M}}\hyperlink{ompreceq}{\preceq}v_{\mathbf{N}}$. For the corresponding result in the unit disk case we refer to Theorem \ref{Thm35disk}. Here, the crucial growth relation \hyperlink{ompreceq}{$\preceq$} between weight functions is defined as follows: Let $v^1$ and $v^2$ be given, then write $v^1\hypertarget{ompreceq}{\preceq}v^2$ if
\begin{equation}\label{bigOrelation}
v^2(t)=O(v^1(t))\;\text{as}\;t\rightarrow+\infty;	
\end{equation}
i.e. $v^2(t)\le Cv^1(t)$ for some $C\ge 1$ and all $t\in[0,+\infty)$ (since both weights are strictly positive). Moreover, write that $v^1$ and $v^2$ are \emph{equivalent,} denoted by $v^1\hypertarget{sim}{\sim}v^2$, if
$$v^1\hyperlink{ompreceq}{\preceq}v^2\;\text{and}\;v^2\hyperlink{ompreceq}{\preceq}v^1.$$
Relation \hyperlink{sim}{$\sim$} is abbreviated by the symbol $\approx$ in \cite{Bonet2022survey} and by $\asymp$ in \cite{AbakumovDoubtsov15} and \cite{AbakumovDoubtsov18}. Similarly, in the unit disk case we write $v^1_{\mathbb{D}}\hyperlink{ompreceq}{\preceq}v^2_{\mathbb{D}}$ if $v^2_{\mathbb{D}}(t)\le Cv^1_{\mathbb{D}}(t)$ for some $C\ge 1$ and all $t\in[0,1)$.

%\begin{lemma}\label{usefulrootlemma}
%Let $\mathbf{M}\in\RR_{>0}^{\NN}$ be given.
%\begin{itemize}
%\item[$(i)$] If $\mathbf{M}$ is log-convex, then for any $k\in\NN$ the map $p\mapsto\left(\frac{M_p}{M_k}\right)^{\frac{1}{p-k}}$ is non-decreasing on $(k,+\infty)$.

%\item[$(ii)$] If $\lim_{p\rightarrow+\infty}(M_p)^{1/p}=+\infty$, then
%$$\forall\;k\in\NN:\;\;\;\lim_{p\rightarrow+\infty}\left(\frac{M_p}{M_k}\right)^{\frac{1}{p-k}}=+\infty.$$
%\end{itemize}
%\end{lemma}

%\demo{Proof}
%$(i)$ For any $k<p$ we have
%$$\left(\frac{M_{p+1}}{M_k}\right)^{\frac{1}{p+1-k}}\ge\left(\frac{M_p}{M_k}\right)^{\frac{1}{p-k}}\Leftrightarrow\frac{M_{p+1}}{M_k}\ge\frac{M_p}{M_k}\left(\frac{M_p}{M_k}\right)^{\frac{1}{p-k}}\Leftrightarrow\mu_{p+1}^{p-k}\ge\frac{M_p}{M_k}=\mu_{k+1}\cdots\mu_p,$$
%and which is clear since by log-convexity $p\mapsto\mu_p$ is non-decreasing.\vspace{6pt}

%$(ii)$ Let $k\in\NN$ and $h>0$ (large) be arbitrary and fixed, then
%$$\left(\frac{M_p}{M_k}\right)^{\frac{1}{p-k}}\ge h\Leftrightarrow M_p\ge h^{p-k}M_k\Leftrightarrow(M_p)^{1/p}\ge h\left(\frac{M_k}{h^k}\right)^{1/p},$$
%and which holds by assumption and since the right-hand side tends to $h$ as $p\rightarrow+\infty$.
%\qed\enddemo

Finally, let us recall the following useful estimates; see e.g. \cite[$(4.1)$, p. 16]{conjugateweightfunction} and we use $0^0:=1$:
\begin{equation}\label{stirling}
\forall\;p\in\NN:\;\;\;p!\le p^p\le e^pp!.
\end{equation}

\subsection{Associated weight function}\label{assofunctionsection}
For the following definition we refer to \cite[Chapitre I]{mandelbrojtbook}, \cite[Sect. 3]{Komatsu73} and the more recent work \cite{regularnew} in which additional properties and explanations are summarized and non-standard cases are treated. Let $\mathbf{M}=(M_p)_{p\in\NN}\in\RR_{>0}^{\NN}$ be given, then the \emph{associated (weight) function} $\omega_{\mathbf{M}}: [0,+\infty)\rightarrow[0,+\infty)\cup\{+\infty\}$ is defined as follows:
\begin{equation}\label{assofunc}
\omega_{\mathbf{M}}(t):=\sup_{p\in\NN}\log\left(\frac{M_0t^p}{M_p}\right),\qquad t\ge 0,
\end{equation}
with the conventions $0^0:=1$ and $\log(0)=-\infty$. This ensures $\omega_{\mathbf{M}}(0)=0$ and $\omega_{\mathbf{M}}(t)\ge 0$ for any $t\ge 0$ since $\frac{t^0M_0}{M_0}=1$ for all $t\ge 0$. \eqref{assofunc} corresponds to \cite[$(3.1)$]{Komatsu73} and we immediately have that $\omega_{\mathbf{M}}$ is non-decreasing and satisfying $\lim_{t\rightarrow+\infty}\omega_{\mathbf{M}}(t)=+\infty$. However, the fact $\omega_{\mathbf{M}}: [0,+\infty)\rightarrow[0,+\infty)$ (well definedness) is equivalent to assume $\lim_{p\rightarrow+\infty}(M_p)^{1/p}=+\infty$; see \cite[Lemma 2.2]{regularnew} and compare this with the above Definition \ref{defweightsequ}. Now let us recall \cite[Lemma 2.4]{regularnew}.

\begin{lemma}\label{lemma1}
Let $\mathbf{M}=(M_p)_{p\in\NN}$ be a log-convex weight sequence. Then
\begin{equation}\label{lemma1equ}
\omega_{\mathbf{M}}(t)=0,\;\text{for}\;t\in[0,\mu_1],\;\;\;\;\omega_{\mathbf{M}}(t)=\log\left(\frac{M_0t^p}{M_p}\right)\;\text{for}\;t\in[\mu_p,\mu_{p+1}],\;p\ge 1.
\end{equation}
\end{lemma}

Consequently, $\omega_{\mathbf{M}}$ is continuous, non-decreasing and tending faster to infinity than any function $t\mapsto\log(t^p)$; thus $\log(t)=o(\omega_{\mathbf{M}}(t))$ as $t\rightarrow+\infty$. Indeed, in \cite[Lemma 2.4]{regularnew} one should read $\log(t^p)$ instead of $\log(t)^p$ and we give a short proof (for $\mathbf{M}$ being a log-convex weight sequence and by taking into account \eqref{lemma1equ}):

$\log(t)=o(\omega_{\mathbf{M}}(t))$ precisely means that for all $p\in\NN_{>0}$ one has $t^p\le\exp(\omega_{\mathbf{M}}(t))$ for all $t$ large enough. Let $p\in\NN_{>0}$ be fixed, then for all $t\in[\mu_{\ell},\mu_{\ell+1}]$ (with $\ell>p$) it is required to get $t^p\le\frac{M_0t^{\ell}}{M_{\ell}}\Leftrightarrow\frac{M_{\ell}}{M_0}\le t^{\ell-p}$ and for this estimate it suffices to have (with $t=\mu_{\ell}$) that $\frac{M_{\ell}}{M_0}=\mu_1\cdots\mu_{\ell}\le\mu_{\ell}^{\ell-p}\Leftrightarrow\mu_1\cdots\mu_p\le\prod_{j=p+1}^{\ell}\frac{\mu_{\ell}}{\mu_j}$. By log-convexity, since all factors in the product on the right-hand side are bounded by below by $1$, one gets $\prod_{j=p+1}^{\ell}\frac{\mu_{\ell}}{\mu_j}\ge\frac{\mu_{\ell}}{\mu_{p+1}}$ and so it suffices to have $\mu_1\cdots\mu_{p+1}=\frac{M_{p+1}}{M_0}\le\mu_{\ell}$. This estimate is clear for all $\ell$ sufficiently large, equivalently for all $t$ large enough, since $\mu_{\ell}\rightarrow+\infty$ as $\ell\rightarrow+\infty$. (If $\mu_{\ell}=\mu_{\ell+1}$, then necessarily $t=\mu_{\ell}$ and hence $\frac{M_0t^{\ell}}{M_{\ell}}=\frac{M_0t^{\ell+1}}{M_{\ell+1}}$.)\vspace{6pt}

Log-convexity is not restricting the generality of this argument: If $\mathbf{M}$ is a weight sequence (according to Definition \ref{defweightsequ}), then the geometric regularizing procedure introducing $\omega_{\mathbf{M}}$ gives, in particular, that $\omega_{\mathbf{M}}=\omega_{\mathbf{M}^{\on{lc}}}$ with $\mathbf{M}^{\on{lc}}$ denoting the \emph{log-convex minorant} of $\mathbf{M}$ and $\mathbf{M}^{\on{lc}}$ is a log-convex weight sequence; see \cite[Chapitre I]{mandelbrojtbook} and \cite[Sect. 3.]{regularnew}. Finally, if $\mathbf{M}$ satisfies
$$\mathbf{M}_{\iota}:=\liminf_{p\rightarrow+\infty}(M_p/M_0)^{1/p}<+\infty,$$
then $\omega_{\mathbf{M}}(t)=+\infty$ for all $t>\mathbf{M}_{\iota}$ (see \cite[Lemma 2.2, $(2.3)$]{regularnew}) which makes $\log(t)=o(\omega_{\mathbf{M}}(t))$ formally trivial.\vspace{6pt}

Altogether this gives the fact that for any (log-convex) weight sequence $\mathbf{M}$ the function $v_{\mathbf{M}}$ is a (radial) weight according to the basic definition in the introduction (rapidly decreasing amounts to $\log(t)=o(\omega_{\mathbf{M}}(t))$).\vspace{6pt}

Let $\mathbf{M}$ be a log-convex weight sequence, then one can compute $\mathbf{M}$ by involving $\omega_{\mathbf{M}}$ as follows, see \cite[Chapitre I, 1.4, 1.8]{mandelbrojtbook}, \cite[Prop. 3.2]{Komatsu73} and more recently \cite[Sect. 3]{regularnew}:
\begin{equation}\label{Prop32Komatsu}
M_p=\sup_{t\ge 0}\frac{t^p}{\exp(\omega_{\mathbf{M}}(t))}=\sup_{t\ge 0}t^pv_{\mathbf{M}}(t),\;\;\;p\in\NN.
\end{equation}
If $\mathbf{M}$ is a weight sequence which is not log-convex, then the right-hand side of \eqref{Prop32Komatsu} yields $\mathbf{M}^{\on{lc}}_p$ and then $\mathbf{M}^{\on{lc}}$ is a log-convex weight sequence.

Moreover, for any log-convex weight sequence let us introduce the \emph{counting function}
\begin{equation*}\label{counting}
\Sigma_{\mathbf{M}}(t):=|\{p\ge 1:\;\;\;\mu_p\le t\}|,
\end{equation*}
and recall the following integral representation formula
\begin{equation}\label{assointrepr}
\forall\;t\ge 0:\;\;\;\omega_{\mathbf{M}}(t)=\int_0^t\frac{\Sigma_{\mathbf{M}}(u)}{u}du=\int_{\mu_1}^t\frac{\Sigma_{\mathbf{M}}(u)}{u}du;
\end{equation}
see \cite[1.8. III]{mandelbrojtbook}, \cite[$(3.11)$]{Komatsu73} and \cite[Lemma 2.5]{regularnew}. Consequently, \begin{equation}\label{assovanishing}
\forall\;t\in[0,\mu_1]:\;\;\;\omega_{\mathbf{M}}(t)=0,
\end{equation}
and so for any $\mathbf{M}\in\hyperlink{LCset}{\mathcal{LC}}$ the function $\omega_{\mathbf{M}}$ vanishes on $[0,1]$ (see also \eqref{lemma1equ}). Thus in this case $v_{\mathbf{M},\mathbb{D}}(0)=1$ holds; hence both $v_{\mathbf{M}}$ and $v_{\mathbf{M},\mathbb{D}}$ are even normalized.

Moreover, in the unit disk case, by definition \eqref{weightfctunitdisk} we get for all $t\in[0,1)$
\begin{equation}\label{weightfctunitdisk1}
v_{\mathbf{M},\mathbb{D}}(t)=\exp\left(\inf_{p\in\NN}-\log\left(\frac{1/(1-t)^p}{M_p}\right)\right)=\exp\left(\inf_{p\in\NN}\log\left(\frac{M_p}{1/(1-t)^p}\right)\right)=\inf_{p\in\NN}M_p(1-t)^p,
\end{equation}
and this identity should be compared with the known and crucial auxiliary function $h_{\mathbf{M}}:t\mapsto\inf_{p\in\NN}t^pM_p$ appearing in the literature dealing with ultradifferentiable and ultraholomorphic functions; see e.g. \cite{almostanalytic} and \cite{optimalflat23}.\vspace{6pt}

For concrete applications, i.e. when one is (only) interested in the associated function w.l.o.g. one can assume that the defining sequence is log-convex since, as mentioned before, one has $\omega_{\mathbf{M}}=\omega_{\mathbf{M}^{\on{lc}}}$. This comment applies to the whole setting in this paper; see also \cite[Rem. 2.5 \& 2.6]{weightedentireinclusion1}. Moreover, recall that $\mathbf{M}\in\hyperlink{LCset}{\mathcal{LC}}$ has been the standard assumption in \cite{solidassociatedweight} and in view of \cite[Rem. 2.5 \& 2.6]{weightedentireinclusion1} for any log-convex weight sequence $\mathbf{M}$ one can find $\mathbf{N}\in\hyperlink{LCset}{\mathcal{LC}}$ such that $H^{\infty}_{v_{\mathbf{M}}}(\CC)=H^{\infty}_{v_{\mathbf{N}}}(\CC)$ as Banach spaces. Similarly, by \cite[Rem. 2.5 \& 2.6]{weightedentireinclusion1} this equality holds for $H^{0}_{v_{\mathbf{M}}}(\CC)$ and the weighted analytic disk spaces $H^{\infty}_{v_{\mathbf{M},\mathbb{D}}}(\DD)$, $H^{0}_{v_{\mathbf{M},\mathbb{D}}}(\DD)$. Summarizing, \cite[Rem. 2.5 \& 2.6]{weightedentireinclusion1} imply the fact that the class \hyperlink{LCset}{$\mathcal{LC}$} is not restricting the generality of the weighted analytic setting. In view of this observation we point out:

\begin{remark}\label{normrem}
For any weight function $v:[0,+\infty)\rightarrow(0,+\infty)$ in the weighted entire setting there exists a normalized $v^n$ such that $v^n(t)=v(t)$ for all $t$ sufficiently large. Indeed, fix $t_0>1$ such that $v(t_0)<\min\{1,v(0)\}(<t_0)$ and such $t_0$ exists since $v(t)$ tends to $0$ as $t\rightarrow+\infty$. Then put $v^n(t):=1$ for $t\in[0,1]$, $v^n(t):=v(t)$ for all $t\ge t_0$ and, finally $v^n(t):=\frac{v(t_0)-1}{t_0-1}t+\frac{t_0-v(t_0)}{t_0-1}$ for all $t\in(1,t_0)$; i.e. affine linear interpolation. $v^n$ satisfies all required conditions since the slope $\frac{v(t_0)-1}{t_0-1}$ is negative.

Similarly, for $v_{\mathbb{D}}$ fix $t_0\in(0,1)$ such that $v_{\mathbb{D}}(t_0)<\min\{1,v_{\mathbb{D}}(0)\}$, put $v^n_{\mathbb{D}}(0):=1$, $v^n_{\mathbb{D}}(t):=v_{\mathbb{D}}(t)$ for all $t\ge t_0$ and $v_{\mathbb{D}}^n(t):=\frac{v_{\mathbb{D}}(t_0)-1}{t_0}t+1$ for all $t\in(0,t_0)$.
\end{remark}

We finish this section with the following technical observation:

\begin{lemma}\label{MG1equivlemma}
Let $\mathbf{M}$ be a log-convex weight sequence. Then the following are equivalent:
\begin{itemize}
\item[$(i)$] $\mathbf{M}\hyperlink{approx}{\approx}\mathbf{G}^1$ holds.

\item[$(ii)$] The sequence of quotients $\mu=(\mu_p)_p$ satisfies
\begin{equation}\label{MG1equivlemmaequ}
0<\inf_{p\in\NN_{>0}}\frac{\mu_p}{p}\le\sup_{p\in\NN_{>0}}\frac{\mu_p}{p}<+\infty.
\end{equation}

\item[$(iii)$] It holds that
$$\exists\;C\ge 1\;\forall\;t\ge 0:\;\;\;\omega_{\mathbf{G}^1}(t/C)\le\omega_{\mathbf{M}}(t)\le\omega_{\mathbf{G}^1}(Ct).$$

\item[$(iv)$] It holds that
$$\exists\;C\ge 1\;\forall\;t\ge 0:\;\;\;v_{\mathbf{G}^1}(Ct)\le v_{\mathbf{M}}(t)\le v_{\mathbf{G}^1}(t/C).$$

\item[$(v)$] It holds that
$$\exists\;C,C_1\ge 1\;\forall\;t\ge 0:\;\;\;\frac{1}{C_1}\exp(-Ct)\le v_{\mathbf{M}}(t)\le C_1\exp(-t/C).$$
\end{itemize}
\end{lemma}

\demo{Proof}
$(i)\Rightarrow(ii)$ Since $M_0=1=G^1_0$ equivalence means $C^{-p}p!\le M_p\le C^pp!$ for some $C\ge 1$ and all $p\in\NN$. Then note that $\mathbf{G}^1$ obviously satisfies \hyperlink{mg}{$(\on{mg})$}, and hence $\mathbf{M}$ too since this property is clearly preserved under equivalence. Recall that for any log-convex weight sequence $\mathbf{M}$ condition \hyperlink{mg}{$(\on{mg})$} is equivalent to \eqref{quotientrootcompprop} and this implies $\mu_p\le A(M_p)^{1/p}\le ACp!^{1/p}\le ACp$ for all $p\in\NN_{>0}$ and hence the second estimate in \eqref{MG1equivlemmaequ}. Moreover, \eqref{stirling}, log-convexity and $M_0=1$ give $p\le ep!^{1/p}\le eC(M_p)^{1/p}\le eC\mu_p$ for all $p\in\NN_{>0}$ and so the first estimate in \eqref{MG1equivlemmaequ}.\vspace{6pt}

$(ii)\Rightarrow(i)$ We have $C^{-1}\le\frac{\mu_p}{p}\le C$, equivalently $pC^{-1}\le\mu_p\le pC$, for some $C\ge 1$ and all $p\in\NN_{>0}$. Since $M_p=\prod_{j=1}^p\mu_j$ we get $p!C^{-p}\le M_p\le p!C^p$ for all $p\in\NN$ and this is precisely $\mathbf{M}\hyperlink{approx}{\approx}\mathbf{G}^1$.\vspace{6pt}

$(i)\Rightarrow(iii)$ $\mathbf{M}\hyperlink{approx}{\approx}\mathbf{G}^1$ implies $\frac{(t/C)^p}{M_p}\le\frac{t^p}{p!}\le\frac{(tC)^p}{M_p}$ for some $C\ge 1$, and for all $p\in\NN$ and all $t\ge 0$. Then, by definition of associated functions, we obtain the desired relation.\vspace{6pt}

$(iii)\Rightarrow(i)$ We apply \eqref{Prop32Komatsu} and get for all $p\in\NN$:
$$M_p=\sup_{t\ge 0}\frac{t^p}{\exp(\omega_{\mathbf{M}}(t))}\le\sup_{t\ge 0}\frac{t^p}{\exp(\omega_{\mathbf{G}^1}(C^{-1}t))}=C^p\sup_{s\ge 0}\frac{s^p}{\exp(\omega_{\mathbf{G}^1}(s))}=C^pG^1_p,$$
and the second (converse) estimate follows similarly.\vspace{6pt}

$(iii)\Leftrightarrow(iv)$ follows immediately by definition \eqref{weights}.\vspace{6pt}

$(iv)\Leftrightarrow(v)$ holds by the known equivalence $\omega_{\mathbf{G}^1}\hyperlink{sim}{\sim}\id: t\mapsto t$; for convenience of the reader we repeat the argument and refer also to the case $s=1$ in \cite[Ex. 2.9]{genLegendreconj}: First, in view of Lemma \ref{lemma1} it holds that $\omega_{\mathbf{G}^1}(t)=0$ for all $0\le t\le 1$ and when $p\le t\le p+1$ for some $p\in\NN_{>0}$, then $\omega_{\mathbf{G}^1}(t)=\log\left(\frac{t^p}{p!}\right)$. This implies for all $p\le t\le p+1$ that
$$\log\left(\frac{p^p}{p!}\right)\le\omega_{\mathbf{G}^1}(t)\le\log\left(\frac{(p+1)^{p}}{p!}\right).$$
$\log\left(\frac{p^p}{p!}\right)\ge\frac{p+1}{C}\Leftrightarrow p^p\ge p!\exp(\frac{p+1}{C})$ is valid for some $C\ge 1$ and all $p\in\NN_{>0}$ by \emph{Stirling's formula} and, similarly, $\log\left(\frac{(p+1)^p}{p!}\right)\le pD\Leftrightarrow (p+1)^p\le p!\exp(pD)$ for some $D\ge 1$ and all $p\in\NN_{>0}$. Finally, note that $\frac{p+1}{C}\ge\frac{t}{C}$ and $pD\le tD$ for all such $t$.
\qed\enddemo

\subsection{Essential weight functions}\label{essentialsection}
We recall the notions from \cite[Def. 1.1, Sect. 1.B]{BierstedtBonetTaskinen98}; see also \cite[Sect. 3]{Bonet2022survey} and the comments and statements in \cite[Sect. 3.1]{weightedentireinclusion1}. Let a weight $v$ be given, then consider the function (radial growth condition) $w_v: \CC\rightarrow(0,+\infty)$ given by
\begin{equation}\label{corrgrowthcond}
w_v:=\frac{1}{v}
\end{equation}
and set
$$B_{w_v}(\CC):=\{f\in H(\CC): |f(z)|\le w_v(|z|),\;\forall\;z\in\CC\}.$$
The function $\widetilde{w_v}:\CC\rightarrow[0,+\infty)$ \emph{associated with} $w_v$ is defined by
\begin{equation}\label{assoweightw}
\widetilde{w_v}(z):=\sup\{|f(z)|: f\in B_{w_v}(\CC)\},\;\;\;z\in\CC,
\end{equation}
and \cite[Observation 1.5]{BierstedtBonetTaskinen98} gives that $\widetilde{w_v}$ is again radial. Finally, set
\begin{equation}\label{assoweightv}
\widetilde{v}:=\frac{1}{\widetilde{w_v}}.
\end{equation}

We recall the following definition; see e.g. \cite[Sect. 3]{Bonet2022survey}.
\begin{definition}\label{essentialdef}
A weight $v$ is called \emph{essential} if $v\hyperlink{sim}{\sim}\widetilde{v}$.
\end{definition}
Similarly these notions can be introduced for arbitrary (open and connected) sets $G$ and even when using a more general setting for $v$ being a weight; see e.g. \cite{BierstedtBonetTaskinen98}, \cite{BonetDomanskiLindstroemTaskinen98} (there $v$ is only assumed to be continuous and strictly positive). In \cite[Prop. 1.3]{BonetDomanskiLindstroemTaskinen98} it has been shown that isometrically $H^{\infty}_{v}(\mathbb{D})=H^{\infty}_{\widetilde{v}}(\mathbb{D})$ but the proof there gives $H^{\infty}_{v}(\CC)=H^{\infty}_{\widetilde{v}}(\CC)$ too; see also \cite[Observ. 1.12]{BierstedtBonetTaskinen98}.\vspace{6pt}

Let $\mathbf{M}$ be a weight sequence and accordingly we write $\widetilde{v_{\mathbf{M}}}$. For any log-convex weight sequence $\mathbf{M}$ the weight $v_{\mathbf{M}}$ is essential for the weighted entire case and even $v_{\mathbf{M}}=\widetilde{v_{\mathbf{M}}}$ holds; see \cite[Lemma 2.16, Thm. 3.4 $(i)$]{weightedentireinclusion1}.

\subsection{(Log-)convex weight functions}\label{convexsection}
For the following notions we also refer to \cite[Sect. 2.6]{weightedentireinclusion1}.

\begin{definition}\label{admissdef}
A weight $v$ is called \emph{convex} if
\begin{equation}\label{vconvexity}
t\mapsto-\log(v(e^t))\;\;\;\text{is convex on}\;\RR.
\end{equation}
Analogously, convexity is defined for $v_{\mathbb{D}}:[0,1)\rightarrow(0,+\infty)$; see e.g. \cite[p. 5, $(4)$]{Bonet2022survey}.
\end{definition}

Note that in \cite{Bonet2022survey}, and occasionally at several other places in the literature, a weight with \eqref{vconvexity} is called \emph{log-convex.} In the weighted entire setting we set
\begin{equation}\label{omegafromv}
\omega^v(t):=-\log(v(t)),\;\;\;t\in[0,+\infty),
\end{equation}
and $\varphi_{\omega^v}:=\omega^v\circ\exp$; see \cite[$(2.19)$ \& $(2.20)$]{weightedentireinclusion1}. Analogously, one introduces $\omega^{v_{\mathbb{D}}}$ and $\varphi_{\omega^{v_{\mathbb{D}}}}$ and so (log-)convexity for $v$ resp. $v_{\mathbb{D}}$ does precisely mean that $\varphi_{\omega^v}$ resp. $\varphi_{\omega^{v_{\mathbb{D}}}}$ is convex. Next we recall that $\omega^v$ satisfies some basic growth properties usually considered for weight functions in the sense of Braun-Meise-Taylor, see \cite{BraunMeiseTaylor90} and \cite[Lemma 2.11]{weightedentireinclusion1}. Indeed, $v$ is (log-)convex if and only if $\omega^v$ satisfies $(\omega_4)$ in \cite{dissertation} and $v$ is normalized if and only if $\omega^v(t)=0$ for $t\in[0,1]$ resp. if and only if $\varphi_{\omega^v}(t)=0$ for all $t\in(-\infty,0]$ which corresponds to the usual normalization condition in the BMT-setting (see e.g. \cite{BraunMeiseTaylor90} and \cite{compositionpaper}). Moreover,
\begin{itemize}
\item[$(*)$] $v$ is rapidly decreasing if and only if

\item[$(*)$] $\log(t)=o(\omega^v(t))$ as $t\rightarrow+\infty$ (i.e. $(\omega_3)$ in \cite{dissertation}) if and only if

\item[$(*)$] $\lim_{t\rightarrow+\infty}\frac{\varphi_{\omega^v}(t)}{t}=+\infty$.
\end{itemize}
Especially, both $(\omega_3)$ and $(\omega_4)$ are standard in the BMT-setting and we refer again to \cite[Lemma 2.11]{weightedentireinclusion1} and the precise study of the convexity condition in \cite{ultradifferentiablecomparison}. Given a weight sequence $\mathbf{M}$, then in view of \eqref{weights} one has $\omega^{v_{\mathbf{M}}}=\omega_{\mathbf{M}}$ and so $v_{\mathbf{M}}$ is convex if and only if $\varphi_{\omega_{\mathbf{M}}}$ is convex and which holds by definition of the associated weight function; see also \cite[Sect. 2.6]{weightedentireinclusion1} and Lemma \ref{lemma1}.

We prove the analogous result for the unit disk case (for this recall \eqref{weightfctunitdisk1}).

\begin{lemma}\label{disclogconvexlemma}
Let $\mathbf{M}$ be a weight sequence. Then $v_{\mathbf{M},\mathbb{D}}$ is (log-)convex and hence essential.
\end{lemma}

\demo{Proof}
According to Definition \ref{admissdef} and \cite[Sect. 3, p. 5, Ex. 1 $(4)$]{Bonet2022survey} we have to show that $t\mapsto-\log(v_{M,\mathbb{D}}(e^t))$ is convex on $(-\infty,0)$. By \eqref{weightfctunitdisk1} this map is given by $t\mapsto-\inf_{p\in\NN}p\log(1-e^t)+\log(M_p)$. On the interval $(-\infty,0)$ the map $t\mapsto e^t$ is clearly (strictly) convex and hence $t\mapsto 1-e^t$ is (strictly) concave and such that $1-e^t\in(0,1)$ for all $t\in(-\infty,0)$. The function $\log$ is on $(0,+\infty)$ (strictly) increasing and concave, hence $t\mapsto\log(1-e^t)$ is concave as well. For the first step note that if $f$ is concave and $g$ is concave and non-decreasing, then on the interval for which the composition $g\circ f$ is defined one has that $g\circ f$ is concave as well because $(g\circ f)(tx+(1-t)y)\ge g(tf(x)+(1-t)f(y))\ge t(g\circ f)(x)+(1-t)(g\circ f)(y)$ for all $x,y$ and all $t\in[0,1]$. Consequently, $t\mapsto p\log(1-e^t)+\log(M_p)$ is concave too for each $p\in\NN$ and so, finally since the infimum of concave functions is concave (and since the infimum exists), we have shown that $t\mapsto-\inf_{p\in\NN}p\log(1-e^t)+\log(M_p)$ is convex as desired.\vspace{6pt}

And (log-)convexity implies the fact that $v_{\mathbf{M},\mathbb{D}}$ is essential since in the unit disk case a weight function $v$ is essential if and only if it is equivalent to a (log-)convex weight function; see \cite[Thm. 2, p. 5]{Bonet2022survey}, \cite{AbakumovDoubtsov18} and also \cite[Rem. 3.6 $(B)$]{weightedentireinclusion1}.
\qed\enddemo

Using this result we can prove the analogues of \cite[Thm. 3.5]{weightedentireinclusion1} and \cite[Thm. 3.2]{weightedentireinclusion2} for the unit disk case as well:

\begin{theorem}\label{Thm35disk}
Let $\mathbf{M}$ and $\mathbf{N}$ be log-convex weight sequences. Then the following are equivalent:
\begin{itemize}
\item[$(i)$] The sequences satisfy \eqref{equ14}; i.e.
$$\exists\;A\ge 1\;\forall\;p\in\NN:\;\;\;N_p\le AM_p.$$
\item[$(ii)$] The weights satisfy
$v_{\mathbf{M},\mathbb{D}}\hyperlink{ompreceq}{\preceq}v_{\mathbf{N},\mathbb{D}}$.
\item[$(iii)$] It holds that $H^{\infty}_{v_{\mathbf{M},\mathbb{D}}}(\DD)\subseteq H^{\infty}_{v_{\mathbf{N},\mathbb{D}}}(\DD)$ (with continuous inclusion).
\item[$(iv)$] It holds that $H^{0}_{v_{\mathbf{M},\mathbb{D}}}(\DD)\subseteq H^{0}_{v_{\mathbf{N},\mathbb{D}}}(\DD)$ (with continuous inclusion).
\end{itemize}
\end{theorem}

\demo{Proof}
By Lemma \ref{disclogconvexlemma} both weights $v_{\mathbf{M},\mathbb{D}}$ and $v_{\mathbf{N},\mathbb{D}}$ are essential and then $(ii)\Leftrightarrow(iii)$ holds by \cite[Thm. 3.2]{weightedentireinclusion1} (alternatively use \cite[Thm. 3.3]{weightedentireinclusion1} instead), whereas $(ii)\Leftrightarrow(iv)$ is valid by \cite[Thm. 3.1]{weightedentireinclusion2}.\vspace{6pt}

For $(i)\Leftrightarrow(ii)$ we follow \cite[Prop. 2.8]{weightedentireinclusion1}:

$(ii)$ implies via \eqref{weightfctunitdisk} that
$\exp\left(-\omega_{\mathbf{N}}\left(\frac{1}{1-t}\right)\right)\le A\exp\left(-\omega_{\mathbf{M}}\left(\frac{1}{1-t}\right)\right)$ for some $A\ge 1$ and all $t\in[0,1)$ and hence $\exp(-\omega_{\mathbf{N}}(s))\le A\exp(-\omega_{\mathbf{M}}(s))$ for all $s\in[1,+\infty)$. When replacing $A$, if necessary, by $A':=\max\{A,\exp(\omega_{\mathbf{M}}(1))\}$, then the above estimate holds for all $s\ge 0$ with $A'$ and hence \eqref{Prop32Komatsu} implies \eqref{equ14}. Conversely, \eqref{equ14} gives $\omega_{\mathbf{M}}(s)\le\omega_{\mathbf{N}}(s)+\log(A)$ for some $A\ge 1$ and all $s\ge 0$. Then restrict this estimate to all $s\in[1,+\infty)$ and the substitution $t:=1-\frac{1}{s}$ immediately gives by definition $v_{\mathbf{M},\mathbb{D}}\hyperlink{ompreceq}{\preceq}v_{\mathbf{N},\mathbb{D}}$.
\qed\enddemo

We close this section by transferring Remark \ref{normrem} in the weighted entire setting to (log-)convex weights. This procedure requires a more sophisticated construction in order to preserve the crucial convexity condition and the technique emphasizes the basic growth property ``rapidly decreasing'' on $v$. Note that analogously in the weight sequence setting, by \cite[Rem. 2.5 \& 2.6]{weightedentireinclusion1} for any log-convex weight sequence $\mathbf{M}$ we can construct an equivalent sequence $\mathbf{N}\in\hyperlink{LCset}{\mathcal{LC}}$ such that the corresponding sequences of quotients $\mu_p$ resp. $\nu_p$ coincide for all $p$ sufficiently large and, moreover, $\mathbf{M}=\mathbf{N}$ if already $\mathbf{M}\in\hyperlink{LCset}{\mathcal{LC}}$. Consequently, the weighted entire spaces (both with $O$- and $o$-growth restrictions) defined in terms of $v_{\mathbf{M}}$ and $v_{\mathbf{N}}$ coincide and, moreover, by \eqref{weights} and the above comments (recall \eqref{assovanishing}) we have that $v_{\mathbf{N}}$ is a (log-)convex and normalized weight function (cf. also \cite[Lemma 2.11]{weightedentireinclusion1}).

\begin{lemma}\label{normconvexrem}
Let $v:[0,+\infty)\rightarrow(0,+\infty)$ be a (log-)convex weight function. Then there exists a normalized (log-)convex weight function $v^n$ such that $v^n(t)=v(t)$ for all $t$ sufficiently large.
\end{lemma}

Since $v^n(t)=v(t)$ for all large $t$ any (relevant) growth condition is preserved since the crucial information is ``located at infinity'' (i.e. $t\rightarrow+\infty$); this comment holds both for the weighted entire setting when considering $v$ and for the ultradifferentiable BMT-weight function setting when involving $\omega^v$ (cf. \eqref{omegafromv}). Please note that with normalization we mean $v^n(t)=0$ for all $t\in[0,1]$ which does not exclude that $v^n(t)=0$ for all $t\in[0,t_0]$ with $t_0>1$; compare this in the weight sequence setting with \eqref{assovanishing}.

In any case, strict positivity of $v$, and Lemma \ref{normconvexrem} resp. Remark \ref{normrem} imply that $v$ and $v^n$ are related by
\begin{equation}\label{normalizedweightrel}
\exists\;C\ge 1\;\forall\;t\ge 0:\;\;\;C^{-1}v(t)\le v^n(t)\le Cv(t).
\end{equation}

\demo{Proof}
Concerning the explicit construction we deal with $\varphi_{\omega^v}=\omega^v\circ\exp$ and distinguish between several cases. Normalization amounts to have $\varphi_{\omega^v}(t)=0$ for all $t\in(-\infty,0]$. By assumption $\varphi_{\omega^v}:\RR\rightarrow\RR$ is convex and we assume that $v$ is not normalized (otherwise set $v:=v^n$). The plan is to define $\widetilde{\varphi}:\RR\rightarrow\RR$ accordingly and then set $$v^n(t):=\exp(-\widetilde{\varphi}(\log(t))),\;\;\;t>0,\hspace{25pt}v^n(0):=0.$$

First, if $\varphi_{\omega^v}(0)=0$ then set
$$\widetilde{\varphi}(t):=0,\;\;\;t\in(-\infty,0],\hspace{25pt}\widetilde{\varphi}(t):=\varphi_{\omega^v}(t),\;\;\;t\in(0,+\infty).$$

Second, if $\varphi_{\omega^v}(0)<0$, then $\varphi_{\omega^v}(t_0)=0$ for some $t_0>0$ since $\omega^v$ is continuous, non-decreasing and $\lim_{t\rightarrow+\infty}\omega^v(t)=+\infty$ (see \eqref{omegafromv}). Then set
$$\widetilde{\varphi}(t):=0,\;\;\;t\in(-\infty,t_0],\hspace{25pt}\widetilde{\varphi}(t):=\varphi_{\omega^v}(t),\;\;\;t\in(t_0,+\infty).$$

By definition it is immediate that in both cases $\widetilde{\varphi}$ is normalized, continuous, non-decreasing, and coinciding with $\varphi_{\omega^v}$ for all $t$ sufficiently large. Convexity follows since $\widetilde{\varphi}$ is given piecewise by two convex functions, the constant zero function $\theta$ and $\varphi_{\omega^v}$. Formally, we have to check $\widetilde{\varphi}(\lambda x+(1-\lambda)y)\le\lambda\widetilde{\varphi}(x)+(1-\lambda)\widetilde{\varphi}(y)$ for all $x,y\in\RR$ and $\lambda\in[0,1]$ and it remains to consider $x<0<y$ resp. $x<t_0<y$. We focus on the first situation since the second one is analogous replacing $0$ by $t_0$. When $\lambda x+(1-\lambda)y\le 0$, then $\widetilde{\varphi}(\lambda x+(1-\lambda)y)=0\le(1-\lambda)\varphi_{\omega^v}(y)=\lambda\widetilde{\varphi}(x)+(1-\lambda)\widetilde{\varphi}(y)$ is clear and when $0<\lambda x+(1-\lambda)y$, then $\widetilde{\varphi}(\lambda x+(1-\lambda)y)=\varphi_{\omega^v}(\lambda x+(1-\lambda)y)\le\lambda\varphi_{\omega^v}(x)+(1-\lambda)\varphi_{\omega^v}(y)\le(1-\lambda)\varphi_{\omega^v}(y)=\lambda\widetilde{\varphi}(x)+(1-\lambda)\widetilde{\varphi}(y)$ since $\varphi_{\omega^v}(x)\le 0$ (note that $\varphi_{\omega^v}$ is non-decreasing) and $\widetilde{\varphi}(x)=0$.\vspace{6pt}

Finally, the third case, when $\varphi_{\omega^v}(0)>0$, follows by the next Lemma \ref{auxlemma} applied to $\varphi_{\omega^v}:=f$, $\widetilde{\varphi}:=f^n$. Note that in this situation \eqref{auxlemmaequ} precisely means that $v$ is rapidly decreasing, whereas for the two previous cases this property is not required necessarily.
\qed\enddemo

\begin{lemma}\label{auxlemma}
Let $f:\RR\rightarrow\RR$ be given and assume that $f$ is convex, non-decreasing, $f(0)>0$, and
\begin{equation}\label{auxlemmaequ}
\lim_{t\rightarrow+\infty}\frac{f(t)}{t}=+\infty.
\end{equation}
Then there exists a convex function $f^n:\RR\rightarrow\RR$ such that $f^n(t)=0$ for $t\in(-\infty,0]$ and $f^n(t)=f(t)$ for all $t$ sufficiently large.
\end{lemma}

\emph{Note:} Indeed, the point $0$ can be replaced by any $\overline{t}\in\RR$, i.e. $f(\overline{t})>0$ instead of $f(0)>0$, by applying a shift of variables $t\mapsto t-\overline{t}$ and which is not effecting the convexity. Without \eqref{auxlemmaequ} this result becomes wrong. Take, e.g., $f_c:\RR\rightarrow\RR$ with $f_c(t):=c$ for arbitrary $c\in(0,+\infty)$, and fix $t_0\in\RR$, then any function $g_{t_0}$ satisfying $g_{t_0}(t)=0$ for $t\in(-\infty,t_0]$ and $g_{t_0}(t)=f_c(t)$ for all sufficiently large $t$ cannot be convex: Assume that there exists $t_1\in\RR$ (with $t_1>t_0$) such that $g_{t_0}=f_c(t)$ for all $t\ge t_1$, then take $x\in(-\infty,t_0]$, $y\in(t_1,+\infty)$ and $x\lambda+(1-\lambda)y\in(t_1,+\infty)$ which can be achieved if $\lambda\in(0,1)$ is chosen sufficiently close to $0$, then $g_{t_0}(\lambda x+(1-\lambda)y)=f_c(\lambda x+(1-\lambda)y)=c>(1-\lambda)c=(1-\lambda)f_c(y)=\lambda g_{t_0}(x)+(1-\lambda)g_{t_0}(y)$ which is contradicting convexity for $g_{t_0}$.

\demo{Proof}
First, we show that for arbitrary $f:\RR\rightarrow\RR$ condition \eqref{auxlemmaequ} is equivalent to the following one:
\begin{equation}\label{auxlemmaequvar}
\forall\;x\in\RR:\;\;\;\lim_{y\rightarrow+\infty}\frac{f(y)-f(x)}{y-x}=+\infty.
\end{equation}
Fix $x\in\RR$, write $\frac{f(y)-f(x)}{y-x}=\frac{f(y)/y-f(x)/y}{1-x/y}$ and so \eqref{auxlemmaequ} implies \eqref{auxlemmaequvar}. Conversely, \eqref{auxlemmaequvar} implies that for (any) fixed $x\in\RR$ and any $C\ge 1$ there exists $y_C\in\RR$, $y_C>x$, such that $f(y)\ge C(y-x)+f(x)$ for all $y\ge y_C$. This yields $\frac{f(y)}{y}\ge C+\frac{f(x)-Cx}{y}$ for all $y\ge y_C$, thus $\liminf_{y\rightarrow+\infty}\frac{f(y)}{y}\ge C$. Since $C\ge 1$ can be chosen arbitrarily, \eqref{auxlemmaequ} is shown.\vspace{6pt}

Let $x_0<y_0$ be arbitrary, then the straight line connecting the points $(x_0,f(x_0))$ and $(y_0,f(y_0))$ is given by the identity $\ell_{x_0,y_0}: t\mapsto k_0t+d_0$ with $k_0:=\frac{f(y_0)-f(x_0)}{y_0-x_0}\ge 0$ (since $f$ is non-decreasing) and $d_0:=f(x_0)-k_0x_0$. Fix now $x_0>0$, then by \eqref{auxlemmaequvar} we can find $y_0>x_0$ such that $k_0>0$ and then $\ell_{x_0,y_0}(t)=0\Leftrightarrow k_0t+d_0=0\Leftrightarrow t=-\frac{d_0}{k_0}=-\frac{f(x_0)}{k_0}+x_0$. Moreover, in view of \eqref{auxlemmaequvar} we can select $y_0>x_0$ large enough to ensure that $t_0:=-\frac{f(x_0)}{k_0}+x_0\ge 0\Leftrightarrow k_0\ge \frac{f(x_0)}{x_0}$. Also we have $t_0<x_0(<y_0)$ since $f(x_0)\ge f(0)>0$ and so $\frac{f(x_0)}{k_0}>0$. Summarizing, let now $0<x_0<y_0$ be fixed such that $t_0\ge 0$ holds and simply write $\ell$ instead of $\ell_{x_0,y_0}$.

Then define
$$f^n(t):=0,\;\;\;t\in(-\infty,t_0],\hspace{20pt}f^n(t):=\ell(t),\;\;\;t\in(t_0,y_0),\hspace{20pt}f^n(t):=f(t),\;\;\;t\ge y_0,$$
and it remains to verify convexity for $f^n$. Since $f^n$ is defined piecewise in terms of three convex functions, namely the zero function $\theta$, the affine linear function $\ell$ and $f$, the desired estimate is clear for all convex-combinations of points $x,y\in(-\infty,t_0]$, $x,y\in[t_0,y_0]$, $x,y\in[y_0,+\infty)$ and we treat the remaining cases.

\begin{itemize}
\item[$(*)$]First, the case $x\in(-\infty,t_0]$ and $y\in[t_0,y_0]$ is a special case of a step in the proof of Lemma \ref{normconvexrem}.\vspace{6pt}

\item[$(*)$] Second, let $x\in[t_0,y_0]$ and $y\in[y_0,+\infty)$ with $x<y$ and we distinguish:

If $\lambda x+(1-\lambda)y\in[t_0,y_0]$, then $f^n(\lambda x+(1-\lambda)y)=\ell(\lambda x+(1-\lambda)y)\le\lambda xk_0+\lambda d_0+(1-\lambda)f(y)=\lambda f^n(x)+(1-\lambda)f^n(y)$ and the desired inequality is equivalent to $(1-\lambda)(yk_0+d_0)\le(1-\lambda)f(y)\Leftrightarrow yk_0+d_0=\ell(y)\le f(y)$ and this holds by convexity of $f$ (note also that $\lambda\in(0,1)$).

If $\lambda x+(1-\lambda)y\in[y_0,+\infty)$, then by convexity of $f$ one has $f^n(\lambda x+(1-\lambda)y)=f(\lambda x+(1-\lambda)y)\le\lambda f(x)+(1-\lambda)f(y)\le\lambda xk_0+\lambda d_0+(1-\lambda)f(y)=\lambda f^n(x)+(1-\lambda)f^n(y)$ and the desired inequality holds since $f(x)\le xk_0+d_0=\ell(x)$ by convexity of $f$.\vspace{6pt}

\item[$(*)$] Finally, let $x\in(-\infty,t_0]$, $y\in[y_0,+\infty)$, thus $x<y$ and we distinguish:

If $\lambda x+(1-\lambda)y\in(-\infty,t_0]$, then $f^n(\lambda x+(1-\lambda)y)=0\le(1-\lambda)f^n(y)=\lambda f^n(x)+(1-\lambda)f^n(y)$ is trivial.

If $\lambda x+(1-\lambda)y\in[t_0,y_0]$, then $f^n(\lambda x+(1-\lambda)y)=\ell(\lambda x+(1-\lambda)y)=\lambda\ell(x)+(1-\lambda)\ell(y)\le 0+(1-\lambda)\ell(y)\le 0+(1-\lambda)f(y)=\lambda f^n(x)+(1-\lambda)f^n(y)$ and the first estimate holds since $\ell(x)\le 0$ whereas the second one is valid by convexity of $f$.

And if $\lambda x+(1-\lambda)y\in[y_0,+\infty)$, then we can assume $y>y_0$ and get $f^n(\lambda x+(1-\lambda)y)=f(\lambda x+(1-\lambda)y)\le(1-\lambda)f^n(y)=\lambda f^n(x)+(1-\lambda)f^n(y)$ as follows: By convexity of $f$ (on $[y_0,+\infty)$) the point $f(\lambda x+(1-\lambda)y)$ lies below the straight line connecting $(y_0,f(y_0))$ and $(y,f(y))$, having slope $k_1:=\frac{f(y)-f(y_0)}{y-y_0}>0$ since convexity also implies $k_1\ge k_0$. Explicitly, this line is given by $\ell_1: t\mapsto k_1t+d_1$ with $d_1=f(y_0)-k_1y_0$ and having the unique intersection point with the $x$-axis at $t_1=-\frac{d_1}{k_1}$. Since $d_0=f(x_0)-k_0x_0=f(y_0)-k_0y_0\ge d_1$ we have that $x\le t_0=-\frac{d_0}{k_0}\le-\frac{d_1}{k_1}=t_1$ and hence the straight line connecting $(x,0)$ and $(y,f(y))$ lies above the point $f(\lambda x+(1-\lambda)y)$.
\end{itemize}
\qed\enddemo

\subsection{Associated weight sequence}\label{assosequencesection}
The following is taken from \cite[Sect. 6]{solidassociatedweight}; see also \cite[Sect. 2.7]{weightedentireinclusion1}. Let $v$ be a normalized and (log-)convex weight. In \cite[$(6.3)$]{solidassociatedweight} we have introduced and studied the sequence $\mathbf{M}^v=(M^v_p)_{p\in\NN}$ defined by
\begin{equation}\label{vBMTweight1equ1}
M^v_p:=\sup_{t>0}t^pv(t)=\sup_{t>0}\frac{t^p}{\exp(\omega^v(t))},\;\;\;p\in\NN.
\end{equation}
Note that this expression has already been considered on \cite[p. 157]{BierstedtBonetTaskinen98} and, alternatively, one can write
\begin{equation}\label{vBMTweight1equ2}
M^v_p=\|z^p\|_v,\;\;\;p\in\NN.
\end{equation}
Since $v$ is normalized we have $M^v_p=\sup_{t\ge 1}t^pv(t)$. The corresponding sequence of quotients is denoted by $\mu^u=(\mu^u_p)_p$ (and put again $\mu^u_0:=1$).

We are also interested in the weight $t\mapsto\exp(-\omega_{\mathbf{M}^v}(t))(=v_{\mathbf{M}^v}(t))$. This notation is suggested by \eqref{weights} but confusing. Thus when explicitly using this function we denote the abstractly given weight with a different symbol than $v$, say $u$, and write $v_{\mathbf{M}^u}$.

We gather several properties; see \cite[Prop. 6.2]{solidassociatedweight} and \cite[Sect. 2.7]{weightedentireinclusion1}:

\begin{itemize}
\item[$(i)$] $\mathbf{M}^u\in\hyperlink{LCset}{\mathcal{LC}}$,

\item[$(ii)$] $\omega_{\mathbf{M}^u}\hyperlink{sim}{\sim}\omega^u$, more precisely $\exists\;A\ge 1\;\forall\;t\ge 0:$
\begin{equation}\label{omegavequiv}
\frac{1}{A}(v_{\mathbf{M}^u}(t))^{2}=\frac{1}{A}\exp(-2\omega_{\mathbf{M}^u}(t))\le u(t)=\exp(-\omega^u(t))\le\exp(-\omega_{\mathbf{M}^u}(t))=v_{\mathbf{M}^u}(t).
\end{equation}

\item[$(iii)$] If $u=v_{\mathbf{M}}$ with $\mathbf{M}\in\hyperlink{LCset}{\mathcal{LC}}$, and so $\omega^u=\omega_{\mathbf{M}}$, then $\mathbf{M}^u=\mathbf{M}$.
\end{itemize}
However, note that (log-)convexity of $u$ is used in this whole approach only for having the first estimate in \eqref{omegavequiv}, see \cite[Rem. 2.13]{weightedentireinclusion1}, and definition \eqref{vBMTweight1equ1} makes sense for more general functions $u$:

\begin{itemize}
\item[$(a)$] Let $u:[0,+\infty)\rightarrow[0,+\infty)$ be rapidly decreasing, then $M^u_p\in(0,+\infty)$ for any $p\in\NN$ and thus $\mathbf{M}^u$ is well defined. Indeed, this assumption ensures for each $p\in\NN$ the existence of at least one critical point of $f_p: t\mapsto t^pu(t)$, write $t_p$ for it, and hence $f_p$ attains its maximum at $t_p$ which yields $M^u_p=f_p(t_p)=t_p^pu(t_p)$. Note that in general there might exist even infinitely many critical points, see e.g. in the weight sequence setting the proof of \cite[Lemma 4.1]{solidassociatedweight} but, however, since for $M^u_p$ one is only interested in the maximum of $f_p$ we can select any of these possible (critical) points and evaluate $f_p$ on these values. Let us choose a corresponding representative and denote it by $t_p$. Hence, whenever dealing with equalities resp. estimates involving critical points (see e.g. \eqref{solidequ65}) this has to be understood that for each $p\in\NN$ we choose a representative among the set of all possible values of critical points.

    Note: For any $p\in\NN_{>0}$ we can extend \eqref{vBMTweight1equ1} to all $t\ge 0$ since $t=0$ yields the value $0$ for the expression in the supremum; for $p=0$ when using the convention $0^0:=1$ for this extension we require that $u$ is continuous at $0$, whereas the convention $0^0:=0$ yields that continuity at $0$ is superfluous for considering $\sup_{t\ge 0}$.

    Finally, note that for any rapidly decreasing $u$ identity \eqref{vBMTweight1equ2} holds.

\item[$(b)$] Even for arbitrary functions $u:[0,+\infty)\rightarrow[0,+\infty)$ the sequence $\mathbf{M}^u$ is always log-convex:
\begin{align*}
(M^u_p)^2&=(\sup_{t>0}t^pu(t))^2=\sup_{t>0}t^{2p}u^2(t)=\sup_{t>0}t^{p-1}u(t)t^{p+1}u(t)\le\sup_{t>0}t^{p-1}u(t)\sup_{t>0}t^{p+1}u(t)
\\&
=M^u_{p-1}M^u_{p+1}.
\end{align*}

\item[$(c)$] If $u:[0,+\infty)\rightarrow[0,+\infty)$ is continuous and non-increasing, then $M^u_0=\sup_{t>0}u(t)=u(0)$ and so normalization for $\mathbf{M}^u$ fails in general.

\item[$(d)$] By \cite[Rem. 2.14]{weightedentireinclusion1} for any continuous and non-increasing $u:[0,+\infty)\rightarrow[0,+\infty)$ (and such that $u\neq 0$), one has $\sup_{p\in\NN_{>0}}(M^u_p)^{1/p}<+\infty$ if and only if $u(t)=0$ for all sufficiently large $t$. In this case, $u$ is rapidly decreasing but not strictly positive anymore. Indeed, the given arguments in \cite{weightedentireinclusion1}, for which $u(0)=1$ is not necessary, show that $$\sup_{p\in\NN_{>0}}(M^u_p/M^u_0)^{1/p}=\lim_{p\rightarrow+\infty}(M^u_p/M^u_0)^{1/p}=:\mathbf{M}^u_{\iota}\in(0,+\infty)$$
    if and only if $u(t)=0$ for all $t>\mathbf{M}^u_{\iota}$. For this note that $M^u_0=v(0)$ (see $(c)$) and by the log-convexity of $\mathbf{M}^u$ (see $(b)$) the map $p\mapsto(M^u_p/M^u_0)^{1/p}$ is non-decreasing (apply \cite[Lemma 2.0.4]{diploma} to $(M^u_p/M^u_0)_p$). Thus the limit exists in $(0,+\infty]$ and hence $\mathbf{M}^u_{\iota}>0$.

    This observation is consistent in the weight sequence setting by recalling \eqref{weights}, \cite[Lemma 2.2, Sect. 4.2]{regularnew}, \cite[Sect. 2.7 $(iii)$]{weightedentireinclusion1} and the convention $\exp(-\infty)=0$.

    Note that formally the basic assumptions on $u$ allow the non-interesting case $u=0$ and which gives $M^u_p=0$ for all $p\in\NN$, and so $\mathbf{M}^u$ is not a weight sequence. However, even in this case the statement formally holds with the convention $\mathbf{M}^v_{\iota}=0$ and in the weight sequence setting this situation corresponds to the limiting case investigated in \cite[Sect. 4.1]{regularnew}.

\item[$(e)$] Remark \ref{normrem} and Lemma \ref{normconvexrem} have now the following meaning: When $u$ is a given weight function, then by switching to the constructed normalized one $u^n$ we get via \eqref{normalizedweightrel} and the definition in \eqref{vBMTweight1equ1} that
\begin{equation}\label{assocsequnormrelation}
\exists\;C\ge 1\;\forall\;p\in\NN:\;\;\;C^{-1}M^u_p\le M^{u^n}_p\le CM^{u}_p,
\end{equation}
which is even stronger than being equivalent.

\item[$(f)$] Let $u:[0,+\infty)\rightarrow(0,+\infty)$ be a weight function, then the sequence (of representatives) $(t_p)_{p\in\NN}$ and the sequence of quotients $\mu^u$ are related as follows (see \cite[$(6.5)$]{solidassociatedweight}):
\begin{equation}\label{solidequ65}
\forall\;p\in\NN:\;\;\;t_p\le\mu^u_{p+1}\le t_{p+1}.
\end{equation}
We repeat the arguments given in \cite{solidassociatedweight} (see also \cite[$(5.10)$]{orliczpaper}) and see that (log-)convexity and normalization for $u$ are not required for having \eqref{solidequ65}. Indeed, we set $t_0=0$ and for this recall $M^u_0=u(0)$ since $u$ is non-increasing and continuous and hence with $0^0:=1$ we can consider $\sup_{t\ge 0}$ in \eqref{vBMTweight1equ1}; see $(a)$. For $p\in\NN$ we get by definition $(t_{p+1})^{p+1}u(t_{p+1})\ge t^{p+1}u(t)$ for all $t\ge 0$ and then choosing $t:=t_p$ gives $M^u_{p+1}=(t_{p+1})^{p+1}u(t_{p+1})\ge t_pt_p^pu(t_p)=t_pM^u_p$. Conversely, $t_{p+1}t_p^pu(t_p)\ge t_{p+1}t^pu(t)$ for all $t\ge 0$ and choosing $t:=t_{p+1}$ gives $t_{p+1}M^u_p\ge(t_{p+1})^{p+1}u(t_{p+1})=M^u_{p+1}$ (again we used for $p=0$ the convention $0^0:=1$).
\end{itemize}

Summarizing, for any weight function $u$ satisfying $u(0)=1$ the corresponding sequence $\mathbf{M}^u$ is formally a log-convex weight sequence according to Definition \ref{defweightsequ} and $(e)$ should be compared with \cite[Rem. 2.5 \& 2.6]{weightedentireinclusion1} in the weight sequence setting since $\mathbf{M}^{u^n}\in\hyperlink{LCset}{\mathcal{LC}}$. Finally, \eqref{solidequ65} implies that $\lim_{p\rightarrow+\infty}t_p=+\infty$ and the importance of this condition is that the sequence $(t_p)_{p\in\NN}$ is more directly related to the weight $u$ and for concrete examples the value of $t_p$ can be computed explicitly.\vspace{6pt}

We apply Lemma \ref{MG1equivlemma} to $\mathbf{M}^u$ and get the following characterization:

\begin{lemma}\label{MuG1equivlemma}
Let $u$ be a normalized and (log-)convex weight function and let $\mathbf{M}^u$ be the associated weight sequence given by \eqref{vBMTweight1equ1} and $(t_p)_{p\in\NN}$ denotes the sequence of (representatives) of critical points. Then the following are equivalent:
\begin{itemize}
\item[$(i)$] $\mathbf{M}^u\hyperlink{approx}{\approx}\mathbf{G}^1$ is valid.

\item[$(ii)$] It holds that
\begin{equation}\label{MuG1equivlemmaequ}
0<\inf_{p\in\NN_{>0}}\frac{\mu^u_p}{p}\le\sup_{p\in\NN_{>0}}\frac{\mu^u_p}{p}<+\infty.
\end{equation}

\item[$(iii)$] It holds that
\begin{equation}\label{MuG1equivlemmaequ1}
0<\inf_{p\in\NN_{>0}}\frac{t_p}{p}\le\sup_{p\in\NN_{>0}}\frac{t_p}{p}<+\infty.
\end{equation}

\item[$(iv)$] It holds that
$$\exists\;C,C_1\ge 1\;\forall\;t\ge 0:\;\;\;\frac{1}{C_1}\exp(-Ct)\le v_{\mathbf{M}^u}(t)\le C_1\exp(-t/C).$$

\item[$(v)$] It holds that
$$\exists\;C,C_1\ge 1\;\forall\;t\ge 0:\;\;\;\frac{1}{C_1}\exp(-Ct)\le u(t)\le C_1\exp(-t/C).$$
\end{itemize}
Indeed, (log-)convexity is only required for having $(iv)\Leftrightarrow(v)$.
\end{lemma}

\emph{Note:} \cite[Cor. 29 (iii)]{Bonet2022survey} corresponds to $(v)$ in Lemma \ref{MG1equivlemma} for $v_{\mathbf{M}}$, to $(iv)$ in Lemma \ref{MuG1equivlemma} for $v_{\mathbf{M}^u}$ and it is $(v)$ in Lemma \ref{MuG1equivlemma} for the weight $u$.

\demo{Proof}
$(i)\Leftrightarrow(ii)\Leftrightarrow(iv)$ follows directly by Lemma \ref{MG1equivlemma} applied to $\mathbf{M}^u$ and $(ii)\Leftrightarrow(iii)$ holds by \eqref{solidequ65}; for these parts (log-)convexity of $u$ is not necessary.\vspace{6pt}

$(v)\Rightarrow(iv)$ We involve \eqref{omegavequiv} (with constant $A\ge 1$) and have $$v_{\mathbf{M}^u}(t)\le\sqrt{A}(u(t))^{1/2}\le\sqrt{AC_1}\exp(-t/(2C)),\;\;\;v_{\mathbf{M}^u}(t)\ge u(t)\ge\frac{1}{C_1}\exp(-Ct)\ge\frac{1}{C_1}\exp(-2Ct),$$
and so it suffices to choose $C':=2C$ and $C_1':=\max\{C_1,\sqrt{AC_1}\}$.

The converse $(iv)\Rightarrow(v)$ follows similarly:
$$u(t)\ge\frac{1}{A}(v_{\mathbf{M}^u}(t))^2\ge\frac{1}{AC_1^2}\exp(-2Ct),\;\;\;u(t)\le v_{\mathbf{M}^u}(t)\le C_1\exp(-t/C)\le C_1\exp(-t/(2C)),$$
and hence we choose (again) $C':=2C$ and $C_1':=AC_1^2$.
\qed\enddemo

\subsection{Weight function systems}\label{parametersection}
Let $v$ be a weight function and consider
$$v_c: t\mapsto v(ct),\hspace{15pt}v^c: t\mapsto v(t)^c,\hspace{15pt}c>0.$$
Concerning this notation we refer to \cite[Sect. 2.2, $(2.3)$]{weightedentireinclusion1} and in the weight sequence setting we write $v_{\mathbf{M},c}$ resp. $v^c_{\mathbf{M}}$ (cf. \cite[$(2.14)$]{weightedentireinclusion1}). In \cite{weightedentireinclusion1} and \cite{weightedentireinclusion2} we have even investigated weighted spaces of entire functions defined in terms of weight (function) systems w.r.t. the parameter $c>0$ and when being expressed in terms of $v_c$ resp. $v^c$. The system using $v_c$ is called of \emph{dilatation-type} whereas $v^c$ is the \emph{exponential-type} system (see \cite[Sect. 2.2, $(2.4)$ \& $(2.5)$, Def. 2.1; Sect. 2.5]{weightedentireinclusion1}).

Note that $v$ is a weight function if and only if some/any $v_c$ resp. if and only if some/any $v^c$ is so, (log-)convexity is preserved under $v\mapsto v^c$, $v\mapsto^v_c$ for arbitrary $c>0$ whereas normalization is preserved for any $v\mapsto v^c$ but only for $v\mapsto v_c$ with $c\in(0,1]$. Finally, for all $c>0$ and any given (log-convex) weight sequence $\mathbf{M}$ set $\mathbf{M}^c:=(c^pM_p)_{p\in\NN}$ and get (see \cite[Rem. 2.10, $(2.18)$]{weightedentireinclusion1}):
\begin{equation}\label{218equ}
\forall\;c>0\;\forall\;t\ge 0:\;\;\;v_{\mathbf{M},c}(t)=v_{\mathbf{M}}(tc)=v_{\mathbf{M}^{\frac{1}{c}}}(t).
\end{equation}
It holds that $\mathbf{M}^c\hyperlink{approx}{\approx}\mathbf{M}$ for all $c>0$, and $\mathbf{M}^c$ is a (log-convex) weight sequence if and only if $\mathbf{M}$ is so (however, normalization is only clear when $c\ge 1$). The corresponding sequence of quotients $\mu^c$ is given by $\mu^c=(c\mu_p)_{p\in\NN_{>0}}$ (and set $\mu^c_0:=1$). Similarly, for the sequence $\mathbf{m}:=\left(\frac{M_p}{p!}\right)_{p\in\NN}$ we can consider $\mathbf{m}^c$.

\section{The differentiation and integration operator in the weight sequence setting}\label{diffintsection}
We provide a detailed study of properties of the \emph{differentiation operator $D$} and the \emph{integration operator $J$} listed in \cite[Sect. 8]{Bonet2022survey} in terms of the defining weight sequence(s). We distinguish between the weighted entire situation and the unit disk case since some results are different and also between spaces with $O$- and $o$-growth restrictions. As explained in \cite[Sect. 8.1, p. 13]{Bonet2022survey} any (log-)convex function $v$ has a right derivative $v'$ everywhere and so, in particular, this statement is valid for $v_{\mathbf{M}}$ (recall \eqref{weights}). For the explicit computation of $v'_{\mathbf{M}}$ we exploit Lemma \ref{lemma1} (see \eqref{lemma1equ}).

When $\mathbf{M}$ is a weight sequence such that even $\mathbf{m}$ is log-convex, then by \cite[Sect. 2.2, Lemma 2.1]{regularnew} (and \cite[Rem. 4.1]{conjugateweightfunction}) applied to $\mathbf{m}$ we have
\begin{equation}\label{stronglc}
\lim_{p\rightarrow+\infty}\frac{\mu_p}{p}=\lim_{p\rightarrow+\infty}\left(\frac{M_p}{p!}\right)^{1/p}=:\mathbf{m}_{\iota}\in(0,+\infty].
\end{equation}

%Finally, for any $a,d,s>0$ we write $v^{1/s}_{d,a}: t\mapsto\exp(-a(dt)^{1/s})$ and this notation is motivated by \cite[Sect. 5.1]{solidassociatedweight}.

\subsection{The weighted entire case with $O$-growth restriction}\label{diffintsectionOgrowth}
In all considerations in this section it turns out that the relative growth of given $\mathbf{M}$ w.r.t. the sequence $\mathbf{G}^{1}$ is relevant. First, consider the \emph{differentiation operator $D: f\mapsto f'$} and by \cite[Thm. 26]{Bonet2022survey} we have an exhaustive characterization:

\begin{theorem}\label{Diffopthm}
Let $\mathbf{M}$ be a log-convex weight sequence. Then the following are equivalent:
\begin{itemize}
	\item[$(i)$] The differentiation operator $D: H^{\infty}_{v_{\mathbf{M}}}(\CC)\rightarrow H^{\infty}_{v_{\mathbf{M}}}(\CC)$ is continuous.
	
	\item[$(ii)$] $-\log(v_{\mathbf{M}}(t))=O(t)$ as $t\rightarrow+\infty$.
	
	\item[$(iii$)] $\omega_{\mathbf{M}}(t)=O(t)$  as $t\rightarrow+\infty$.
	
	\item[$(iv)$] $\liminf_{p\rightarrow+\infty}\left(\frac{M_p}{p!}\right)^{1/p}>0$ (i.e. $\mathbf{G}^{1}\hyperlink{preceq}{\preceq}\mathbf{M}$).

    \item[$(v)$] $\liminf_{p\rightarrow+\infty}\frac{\mu_p}{p}>0$.
\end{itemize}

Moreover, the following are equivalent:

\begin{itemize}
	\item[$(i)'$] The differentiation operator $D: H^{\infty}_{v_{\mathbf{M}}}(\CC)\rightarrow H^{\infty}_{v_{\mathbf{M}}}(\CC)$ is compact.
	
	\item[$(ii)'$] $-\log(v_{\mathbf{M}}(t))=o(t)$ as $t\rightarrow+\infty$.
	
	\item[$(iii)'$] $\omega_{\mathbf{M}}(t)=o(t)$ as $t\rightarrow+\infty$.
	
	\item[$(iv)'$] $\lim_{p\rightarrow+\infty}\left(\frac{M_p}{p!}\right)^{1/p}=+\infty$ (i.e. $\mathbf{G}^{1}\hyperlink{mtriangle}{\vartriangleleft}\mathbf{M}$).

\item[$(v)'$]  $\lim_{p\rightarrow+\infty}\frac{\mu_p}{p}=+\infty$.
\end{itemize}
\end{theorem}

\emph{Note:} Conditions $(iv)$ and $(iv)'$ are obviously preserved under equivalence of weight sequences and $(iv)$ (resp. $(v)$) is automatically valid if $\mathbf{M}$ is strongly log-convex with $\lim_{p\rightarrow+\infty}\left(\frac{M_p}{p!}\right)^{1/p}=\lim_{p\rightarrow+\infty}\frac{\mu_p}{p}\in(0,+\infty]$ (cf. \eqref{stronglc}).

\demo{Proof}
First, recall that $v_{\mathbf{M}}$ is a radial and (log-)convex weight on $\CC$ according to the requirements in \cite[Sect. 8.1]{Bonet2022survey}. Hence the equivalences $(i)\Leftrightarrow(ii)$ and $(i)'\Leftrightarrow(ii)'$ follow by \cite[Thm. 26]{Bonet2022survey} applied to $v=v_{\mathbf{M}}$, see also \cite[Thm. 2.10]{AbaninThien17} and \cite{AbaninThien18} and note that the symbol for the weight $v$ in the latter works corresponds to $v^{-1}$ in \cite{Bonet2022survey}.

Moreover, by definition \eqref{weights} we immediately get $(ii)\Leftrightarrow(iii)$ and $(ii)'\Leftrightarrow(iii)'$ and therefore it remains to verify $(iii)\Leftrightarrow(iv)\Leftrightarrow(v)$ resp. $(iii)'\Leftrightarrow(iv)'\Leftrightarrow(v)'$.\vspace{6pt}

$(iii)\Rightarrow(iv)$ When $\omega_{\mathbf{M}}(t)=O(t)$ holds, then
$$\exists\;C\ge 1\;\forall\;t\ge 0:\;\;\;\omega_{\mathbf{M}}(t)\le Ct+C,$$
hence \eqref{Prop32Komatsu} implies for all $p\in\NN$:
$$M_p=\sup_{t\ge 0}\frac{t^p}{\exp(\omega_M(t))}\ge e^{-C}\sup_{t\ge 0}\frac{t^p}{\exp(Ct)}=e^{-C}C^{-p}\sup_{s\ge 0}\frac{s^p}{\exp(s)}.$$
Let $p\in\NN_{>0}$ be given and set $f_p(s):=\frac{s^p}{\exp(s)}$, $s\in(0,+\infty)$. Then $f'_p(s)=ps^{p-1}e^{-s}-s^pe^{-s}$ and thus $f'_p(s)=0$ has the unique and global maximum point $s=p$ which gives $\sup_{s\ge 0}\frac{s^p}{\exp(s)}=f_p(p)=\frac{p^p}{e^p}$ for any $p\in\NN_{>0}$ and for $p=0$ clearly $\sup_{s\ge 0}\frac{1}{\exp(s)}=1=\frac{p^p}{e^p}$ by using the convention $0^0:=1$. Consequently, \eqref{stirling} implies
$$\exists\;C\ge 1\;\forall\;p\in\NN:\;\;\;e^{-C}\frac{1}{(eC)^p}\le\frac{M_p}{p^p}\le\frac{M_p}{p!},$$
which gives the conclusion.\vspace{6pt}

$(iv)\Rightarrow(v)$ First, by \eqref{stirling} one gets $\frac{p!^{1/p}}{(p+1)}\ge\frac{p}{e(p+1)}\ge\frac{1}{2e}$, hence $\frac{\mu_{p+1}}{p!^{1/p}}\le 2e\frac{\mu_{p+1}}{p+1}$ for all $p\in\NN_{>0}$. Second, since $\mathbf{M}$ is log-convex and $M_0=1$ we get that $M_p=\mu_1\cdots\mu_p\le\mu_{p+1}^p$ and hence $\frac{(M_p)^{1/p}}{p!^{1/p}}\le\frac{\mu_{p+1}}{p!^{1/p}}$ for all $p\in\NN_{>0}$. By combining both estimates the conclusion follows. Indeed, for the second estimate it suffices to assume the weaker requirement that $p\mapsto(M_p)^{1/p}$ is non-decreasing (and $M_0=1$); see e.g. \cite[Lemma 2.0.4]{diploma}. This is valid by the following equivalences:
\begin{align*}
&(M_p)^{1/p}\le(M_{p+1})^{1/(p+1)}\Leftrightarrow\frac{p!^{1/p}}{(p+1)!^{1/(p+1)}}(M_p/p!)^{1/p}\le(M_{p+1}/(p+1)!)^{1/(p+1)}
\\&
\Leftrightarrow\frac{p!p!^{1/p}}{(p+1)!}\frac{M_p}{p!}(M_p/p!)^{1/p}\le\frac{M_{p+1}}{(p+1)!}\Leftrightarrow(M_p/p!)^{1/p}\le\frac{\mu_{p+1}}{p!^{1/p}}.
\end{align*}

$(v)\Rightarrow(iv)$ By assumption,
$$\exists\;\epsilon>0\;\forall\;p\in\NN_{>0}:\;\;\;\frac{\mu_p}{p}\ge\epsilon\Leftrightarrow M_p\ge\epsilon pM_{p-1}.$$
Thus iteration gives $M_p\ge\epsilon^pp!M_0$ for all $p\in\NN_{>0}$ and for $p=0$ one has the equality $M_0=M_0(=1)$. This implies $(M_p/p!)^{1/p}\ge\epsilon(M_0)^{1/p}$ for all $p\ge 1$ and as $p\rightarrow+\infty$ one infers $\liminf_{p\rightarrow+\infty}(M_p/p!)^{1/p}\ge\epsilon>0$ as desired.\vspace{6pt}

$(iv)\Rightarrow(iii)$ By assumption,
$$\exists\;\epsilon>0\;\forall\;p\in\NN:\;\;\;p!\epsilon^p\le M_p,$$
and note that $p=0$ yields the equality $1=1$. Then $\frac{t^p}{M_p}\le\frac{(t/\epsilon)^p}{p!}$ for all $p\in\NN$ and $t\ge 0$ which gives by definition $\omega_{\mathbf{M}}(t)\le\omega_{\mathbf{G}^1}(t/\epsilon)$ for all $t\ge 0$ (note that $G^1_0=0!=1$). The conclusion follows by $\omega_{\mathbf{G}^1}\hyperlink{sim}{\sim}\id$; recall the proof of $(iv)\Leftrightarrow(v)$ in Lemma \ref{MG1equivlemma}.\vspace{6pt}

$(iii)'\Rightarrow(iv)'$ We even have $\omega_{\mathbf{M}}(t)=o(t)$ and so
$$\forall\;c>0\;\exists\;C\ge 1\;\forall\;t\ge 0:\;\;\;\omega_M(t)\le ct+C,$$
and following the arguments above one infers
$$\forall\;c>0\;\exists\;C\ge 1\;\forall\;p\in\NN:\;\;\;e^{-C}\frac{1}{(ec)^p}\le\frac{M_p}{p^p}\le\frac{M_p}{p!},$$
hence $(iv)'$ is shown.\vspace{6pt}

$(iv)'\Rightarrow(v)'$ holds by the estimate and arguments from $(iv)\Rightarrow(v)$ and for $(v)'\Rightarrow(iv)'$ we proceed as follows: By assumption
$$\forall\;h>0\;\exists\;p_h\in\NN_{>0}\;\forall\;p>p_h:\;\;\;\frac{\mu_p}{p}\ge h.$$
Fix $h>0$ (large), then for all $p>p_h$ we get $$M_p=\mu_1\cdots\mu_{p_h}\cdot\mu_{p_h+1}\cdots\mu_p=M_{p_h}\cdot\mu_{p_h+1}\cdots\mu_p\ge M_{p_h}h^{p-p_h}\frac{p!}{p_h!}$$
and so $(M_p/p!)^{1/p}\ge h\left(\frac{M_{p_h}}{h^{p_h}p_h!}\right)^{1/p}$ for all $p$ sufficiently large, hence $\liminf_{p\rightarrow+\infty}(M_p/p!)^{1/p}\ge h$. Since $h>0$ can be chosen arbitrary large we are done.\vspace{6pt}

$(iv)'\Rightarrow(iii)'$ When $\lim_{p\rightarrow+\infty}\left(\frac{M_p}{p!}\right)^{1/p}=+\infty$, then
$$\forall\;h\ge 1\;\exists\;C\ge 1\;\forall\;p\in\NN:\;\;\;p!h^p\le CM_p,$$
which gives $\omega_{\mathbf{M}}(t)\le\omega_{\mathbf{G}^1}(h^{-1}t)+\log(C)$ for all $t\ge 0$ and so $\omega_{\mathbf{M}}(t)=o(t)$ holds by taking into account again $\omega_{\mathbf{G}^1}\hyperlink{sim}{\sim}\id$.
\qed\enddemo

By involving \eqref{solidequ65}, $(ii)$ in Section \ref{assosequencesection} (see \eqref{omegavequiv}) and the fact that assertions $(iii)$ and $(iii)'$ in Theorem \ref{Diffopthm} are obviously preserved under equivalence \hyperlink{sim}{$\sim$}, we immediately get the following consequence which is extending the characterization from \cite[Thm. 26]{Bonet2022survey} for abstractly given weights $u$ in terms of $\mathbf{M}^u$:

\begin{corollary}\label{Diffopthmcor}
Let $u$ be a normalized and (log-)convex weight function, $\mathbf{M}^u$ the corresponding associated weight sequence and $(t_p)_{p\in\NN}$ denotes the sequence of (representatives of) critical points. Then the following are equivalent:

	\begin{itemize}
		\item[$(i)$] The differentiation operator $D: H^{\infty}_{u}(\CC)\rightarrow H^{\infty}_{u}(\CC)$ is continuous.
		
		\item[$(ii)$] $-\log(u(t))=\omega^u(t)=O(t)$ as $t\rightarrow+\infty$.
		
		\item[$(iii$)] $\omega_{\mathbf{M}^u}(t)=O(t)$ as $t\rightarrow+\infty$.
		
		\item[$(iv)$] $\liminf_{p\rightarrow+\infty}\left(\frac{M^u_p}{p!}\right)^{1/p}>0$ (i.e. $\mathbf{G}^{1}\hyperlink{preceq}{\preceq}\mathbf{M}^u$).

\item[$(v)$] $\liminf_{p\rightarrow+\infty}\frac{\mu^u_p}{p}>0$.

\item[$(vi)$] $\liminf_{p\rightarrow+\infty}\frac{t_p}{p}>0$.

\item[$(vii)$] The differentiation operator $D: H^{\infty}_{v_{\mathbf{M}^u}}(\CC)\rightarrow H^{\infty}_{v_{\mathbf{M}^u}}(\CC)$ is continuous.
	\end{itemize}
	
	Moreover, the following are equivalent:
	
	\begin{itemize}
		\item[$(i)'$] The differentiation operator $D: H^{\infty}_{u}(\CC)\rightarrow H^{\infty}_{u}(\CC)$ is compact.
		
		\item[$(ii)'$] $-\log(u(r))=\omega^u(t)=o(t)$ as $t\rightarrow+\infty$.
		
		\item[$(iii)'$] $\omega_{\mathbf{M}^u}(t)=o(t)$ as $t\rightarrow+\infty$.
		
		\item[$(iv)'$] $\lim_{p\rightarrow+\infty}\left(\frac{M^u_p}{p!}\right)^{1/p}=+\infty$ (i.e. $\mathbf{G}^{1}\hyperlink{mtriangle}{\vartriangleleft}\mathbf{M}^u$).

       \item[$(v)'$] $\lim_{p\rightarrow+\infty}\frac{\mu^u_p}{p}=+\infty$.

       \item[$(vi)'$] $\lim_{p\rightarrow+\infty}\frac{t_p}{p}=+\infty$.

\item[$(vii)'$] The differentiation operator $D: H^{\infty}_{v_{\mathbf{M}^u}}(\CC)\rightarrow H^{\infty}_{v_{\mathbf{M}^u}}(\CC)$ is compact.
	\end{itemize}
\end{corollary}

Using the notation from Section \ref{parametersection} we immediately get the next observation:

\begin{remark}\label{Diffopthmrem}
The following are equivalent:
\begin{itemize}
	\item[$(*)$] $-\log(v(t))=O(t)$ resp. $-\log(v(t))=o(t)$ as $t\rightarrow+\infty$,
	
	\item[$(*)$] $-\log(v^c(t))=O(t)$ resp. $-\log(v^c(t))=o(t)$ as $t\rightarrow+\infty$ for some/any $c>0$,
	
	\item[$(*)$] $-\log(v_c(t))=O(t)$ resp. $-\log(v_c(t))=o(t)$ as $t\rightarrow+\infty$ for some/any $c>0$.
\end{itemize}
Thus the equivalences stated in Theorem \ref{Diffopthm} and Corollary \ref{Diffopthmcor} are preserved under the transformations $v\mapsto v_c$ and $v\mapsto v^c$; recall that (log-)convexity is preserved under both modifications, normalization w.r.t. $v\mapsto v^c$ but $v\mapsto v_c$ only with $c\in(0,1]$. However, in view of Lemma \ref{normconvexrem} in any situation we can switch to a normalized weight and which does not effect the underlying spaces; the corresponding associated weight sequences are then related by \eqref{assocsequnormrelation}.

In view of this comment note that $\mathbf{M}$ is strongly log-convex if and only if some/any $\mathbf{M}^c$ is so.
\end{remark}

We continue with \cite[Thm. 28]{Bonet2022survey} and focus on the \emph{integration operator $J(f)(z):=\int_0^zf(\zeta)d\zeta$, $z\in\CC$:}

\begin{theorem}\label{Intopthm}
Let $\mathbf{M}$ be a log-convex weight sequence. Consider the following assertions:
	\begin{itemize}
		\item[$(i)$] The integration operator $J: H^{\infty}_{v_{\mathbf{M}}}(\CC)\rightarrow H^{\infty}_{v_{\mathbf{M}}}(\CC)$ is continuous.
		
		\item[$(ii)$] $\limsup_{t\rightarrow+\infty}v_{\mathbf{M}}(t)\int_0^t\frac{1}{v_{\mathbf{M}}(r)}dr<+\infty$.
		
		\item[$(iii$)] $\liminf_{t\rightarrow+\infty}v_{\mathbf{M}}(t)\left(\frac{1}{v_{\mathbf{M}}}\right)'(t)>0$.
		
		\item[$(iv)$] $\sup_{p\in\NN_{>0}}\frac{\mu_p}{p}<+\infty$.

\item[$(v)$] $\sup_{p\in\NN_{>0}}\left(\frac{M_p}{p!}\right)^{1/p}<+\infty$ (i.e. $\mathbf{M}\hyperlink{preceq}{\preceq}\mathbf{G}^{1}$).
	\end{itemize}
	
Moreover, consider the following assertions:

	\begin{itemize}
		\item[$(i)'$] The integration operator $J: H^{\infty}_{v_{\mathbf{M}}}(\CC)\rightarrow H^{\infty}_{v_{\mathbf{M}}}(\CC)$ is compact.
		
		\item[$(ii)'$] $\lim_{t\rightarrow+\infty}v_{\mathbf{M}}(t)\int_0^t\frac{1}{v_{\mathbf{M}}(r)}dr=0$.
		
		\item[$(iii)'$] $\lim_{t\rightarrow+\infty}v_{\mathbf{M}}(t)\left(\frac{1}{v_{\mathbf{M}}}\right)'(t)=+\infty$.
		
		\item[$(iv)'$] $\lim_{p\rightarrow+\infty}\frac{\mu_p}{p}=0$.

\item[$(v)'$] $\lim_{p\rightarrow+\infty}\left(\frac{M_p}{p!}\right)^{1/p}=0$ (i.e. $\mathbf{M}\hyperlink{mtriangle}{\vartriangleleft}\mathbf{G}^1$).
	\end{itemize}
Then $(i)\Leftrightarrow(ii)\Leftrightarrow(iii)\Leftrightarrow(iv)\Rightarrow(v)$ and $(i)'\Leftrightarrow(ii)'\Leftrightarrow(iii)'\Leftrightarrow(iv)'\Rightarrow(v)'$ is valid.

Moreover:
\begin{itemize}
\item[$(*)$] When $\mathbf{M}$ satisfies in addition \hyperlink{mg}{$(\on{mg})$}, then the implications $(v)\Rightarrow(iv)$, $(v)'\Rightarrow(iv)'$ are valid and hence equivalence of all assertions.

\item[$(*)$] When in addition even $\left(\frac{M_p}{p!}\right)_{p\in\NN}$ is log-convex, i.e. $\mathbf{M}$ is \emph{strongly log-convex,} then one has $\lim_{p\rightarrow+\infty}\left(\frac{M_p}{p!}\right)^{1/p}=\lim_{p\rightarrow+\infty}\frac{\mu_p}{p}\in(0,+\infty]$ (cf. \eqref{stronglc}). Consequently, on the one hand, this gives directly $(iv)\Leftrightarrow(v)$ and hence again equivalence of the first five assertions, but, on the other hand, both $(iv)'$ and $(v)'$ are violated.
\end{itemize}
\end{theorem}

\emph{Note:}

\begin{itemize}
\item[$(*)$] $(v)$ and $(v)'$ are obviously preserved under equivalence of sequences.

\item[$(*)$] The most prominent example for a strongly log-convex weight sequence is given by $\mathbf{G}^s$, $s>1$.

\item[$(*)$] $(v)'$ is non-standard in the ultrdifferentiable weight sequence setting (cf. \cite{Komatsu73}) but not contradicting the fact that $\mathbf{M}$ is a (log-convex) weight sequence according to Definition \ref{defweightsequ}. A prominent example is provided by the sequences $\mathbf{G}^s$ with $0<s<1$; i.e. consider \emph{small Gevrey indices resp. Gevrey-Markin sequences} and we refer to \cite{microclasses} and \cite{conjugateweightfunction} and the references therein. A different sequence is constructed explicitly in \cite[Prop. 6.9]{conjugateweightfunction}. Moreover, recall that if $p\mapsto\frac{\mu_p}{p}$ is non-increasing, i.e. if $\mathbf{m}$ is \emph{log-concave}, then $\mathbf{M}\hyperlink{mtriangle}{\vartriangleleft}\mathbf{G}^1$ is (also) equivalent to having $\lim_{p\rightarrow+\infty}\frac{\mu_p}{p}=0$; see \cite[Sect. 4.3]{conjugateweightfunction} and \cite[Sect. 2.2]{microclasses} for details. However, in view of \cite[Rem. 4.18]{conjugateweightfunction} this condition already implies \hyperlink{mg}{$(\on{mg})$}.

    As said (log-convex) sequences $\mathbf{M}$ satisfying $\mathbf{M}\hyperlink{mtriangle}{\vartriangleleft}\mathbf{G}^1$ are exotic in the standard setting but nevertheless they play an important role in different contexts, for example when studying the regularity of (weak) solutions of an evolution equation within the Hilbert space setting; we refer again to \cite{microclasses} and \cite{conjugateweightfunction}.
\end{itemize}

\demo{Proof}
First, the equivalences $(i)\Leftrightarrow(ii)\Leftrightarrow(iii)$ and $(i)'\Leftrightarrow(ii)'\Leftrightarrow(iii)'$ hold by applying \cite[Thm. 28]{Bonet2022survey} to $v=v_{\mathbf{M}}$; see also \cite[Thm. 3.8]{AbaninThien17} and compare again the conditions in view of the different notation used in \cite{AbaninThien17}.\vspace{6pt}

We verify $(iii)\Leftrightarrow(iv)$ and $(iii)'\Leftrightarrow(iv)'$. For this note that by Lemma \ref{lemma1} the expression in assertions $(iii)$ and $(iii)'$ takes the following form: Let $t\ge 0$ be such that $\mu_p\le t<\mu_{p+1}$ for some $p\in\NN_{>0}$, then
\begin{align*}
v_{\mathbf{M}}(t)\left(\frac{1}{v_{\mathbf{M}}}\right)'(t)&=\exp(-\omega_{\mathbf{M}}(t))\exp(\omega_{\mathbf{M}}(t))'=\exp(-\omega_{\mathbf{M}}(t))\exp(\omega_{\mathbf{M}}(t))\omega'_{\mathbf{M}}(t)
\\&
=\log\left(\frac{t^p}{M_p}\right)'=\frac{p}{t}.
\end{align*}
$(iii)\Rightarrow(iv)$ By assumption it holds that
\begin{equation}\label{Intopthmequ}
\exists\;\epsilon>0\;\exists\;t_{\epsilon}>0\;\forall\;t\ge t_{\epsilon}:\;\;\;v_{\mathbf{M}}(t)\left(\frac{1}{v_{\mathbf{M}}}\right)'(t)>\epsilon.
\end{equation}
Let $p_0\in\NN_{>0}$ be chosen such that $\mu_{p_0}\ge t_{\epsilon}$ and for this note that $p\mapsto\mu_p$ is non-decreasing and $\lim_{p\rightarrow+\infty}\mu_p=+\infty$ since $\mathbf{M}$ is a log-convex weight sequence. Now let $t\ge\mu_{p_0}$ and so $\mu_p\le t<\mu_{p+1}$ for some $p\ge p_0$. By the above equation we get $\epsilon<\frac{p}{t}\le\frac{p}{\mu_p}$; thus $\frac{\mu_p}{p}<\epsilon^{-1}$ for all $p\ge p_0$ satisfying $\mu_{p+1}>\mu_p$ (and w.l.o.g. one can assume $\mu_{p_0}<\mu_{p_0+1}$). And when $p>p_0$ is such that $\mu_{p-1}<\mu_{p}=\dots=\mu_{p+\ell}<\mu_{p+\ell+1}$ for some $\ell\in\NN_{>0}$, then by assumption and the above equation $\frac{p-1}{t}>\epsilon$ for all $t\in[\mu_{p-1},\mu_p)$. Consequently $\frac{p-1}{\mu_p}\ge\epsilon\Leftrightarrow\epsilon^{-1}\ge\frac{\mu_p}{p-1}$ which implies $\frac{\mu_{p+i}}{p+i}\le\epsilon^{-1}$ for all $i=0,\dots,\ell$ and so the desired estimate is valid for all $p\ge p_0$ when choosing $A:=\epsilon^{-1}$. By enlarging the constant $A\ge 1$ if necessary one gets $\frac{\mu_p}{p}\le A$ for all $p\in\NN_{>0}$ and hence the conclusion follows.\vspace{6pt}

$(iii)'\Rightarrow(iv)'$ If even $\lim_{t\rightarrow+\infty}v_{\mathbf{M}}(t)\left(\frac{1}{v_{\mathbf{M}}}\right)'(t)=+\infty$, then \eqref{Intopthmequ} turns into
\begin{equation*}\label{Intopthmequ1}
\forall\;h>0\;\exists\;t_{h}>0\;\forall\;t\ge t_{h}:\;\;\;v_{\mathbf{M}}(t)\left(\frac{1}{v_{\mathbf{M}}}\right)'(t)>h.
\end{equation*}
Let $h\ge 1$ (large) be fixed and let $p_h\in\NN_{>0}$ be chosen (minimal) such that $\mu_{p_h}\ge t_{h}$ and that $\mu_{p_h}<\mu_{p_h+1}$. Similarly as before one verifies $\frac{\mu_p}{p}\le h^{-1}$ for all $p\ge p_h$. Since $h$ can be chosen arbitrarily large, $\lim_{p\rightarrow+\infty}\frac{\mu_p}{p}=0$ is shown.\vspace{6pt}
		
$(iv)\Rightarrow(iii)$ $\sup_{p\in\NN_{>0}}\frac{\mu_p}{p}<+\infty$ implies $\mu_p\le Dp$ for some $D\ge 1$ and all $p\in\NN_{>0}$. Let $t\ge\mu_1$ be given, then $\mu_p\le t<\mu_{p+1}$ for some $p\in\NN_{>0}$ and therefore, by the above computation, we conclude that $$v_{\mathbf{M}}(t)\left(\frac{1}{v_{\mathbf{M}}}\right)'(t)=\frac{p}{t}\ge\frac{p+1}{2\mu_{p+1}}\ge\frac{1}{2D}>0.$$

$(iv)'\Rightarrow(iii)'$ If even $\lim_{p\rightarrow+\infty}\frac{\mu_p}{p}=0$, then for any $c>0$ we can find $p_c\in\NN_{>0}$ such that $\mu_p\le cp$ for all $p\ge p_c$. Again, w.l.o.g. assume $\mu_{p_c}<\mu_{p_c+1}$ and then, similarly,
$$v_{\mathbf{M}}(t)\left(\frac{1}{v_{\mathbf{M}}}\right)'(t)=\frac{p}{t}\ge\frac{p+1}{2\mu_{p+1}}\ge\frac{1}{2c}$$
holds for all $t\ge\mu_{p_c}$ and when $c\rightarrow 0$ we are done.\vspace{6pt}

$(iv)\Rightarrow(v)$ By \eqref{stirling} and $(M_p/p!)^{1/p}\le\mu_{p+1}/p!^{1/p}$ for all $p\in\NN_{>0}$ (see $(iv)\Rightarrow(v)$ in Theorem \ref{Diffopthm}) we have $(M_p/p!)^{1/p}\le\frac{\mu_{p+1}}{p!^{1/p}}\le e\frac{\mu_{p+1}}{p}\le 2e\frac{\mu_{p+1}}{p+1}$ for all $p\in\NN_{>0}$ which yields $(v)$. And this estimate also gives $(iv)'\Rightarrow(v)'$.\vspace{6pt}

Concerning the supplement $(v)\Rightarrow(iv)$, $(v)'\Rightarrow(iv)'$, first we show that for any log-convex weight sequence $\mathbf{M}$ condition \hyperlink{mg}{$(\on{mg})$} is equivalent to
\begin{equation}\label{quotientrootcomppropgen}
\forall\;s\in\RR\;\exists\;B\ge 1\;\forall\;p\in\NN_{>0}:\;\;\;\mu_pp^s\le B(M_pp!^s)^{1/p}.
\end{equation}
Indeed, \eqref{quotientrootcomppropgen} for $s=0$ precisely corresponds to \eqref{quotientrootcompprop} and hence yields  \hyperlink{mg}{$(\on{mg})$}. Conversely, let $\mathbf{M}$ have \hyperlink{mg}{$(\on{mg})$} and so \eqref{quotientrootcompprop}; i.e. \eqref{quotientrootcomppropgen} for $s=0$. Then distinguish:

If $s>0$, then \eqref{quotientrootcompprop} and \eqref{stirling} give $\mu_pp^s\le A(M_p)^{1/p}p^s=A(M_pp^{ps})^{1/p}\le A(M_p e^{ps}p!^s)^{1/p}=Ae^s(M_pp!^s)^{1/p}$. If $s<0$, then $\mu_pp^s\le A(M_p)^{1/p}p^s=A(M_pp^{ps})^{1/p}\le A(M_pp!^s)^{1/p}$. Applying now \eqref{quotientrootcomppropgen} to $s=-1$ gives the conclusion.
\qed\enddemo

Clearly, $(iv)$ resp. $(iv)'$ holds for some/any $\mathbf{M}^c$ if and only if $\mathbf{M}$ does so and $(v)$ resp. $(v)'$ are preserved under equivalence of sequences. On the other hand, even for abstractly given weight functions $(ii)$ resp. $(ii)'$ are preserved under $v\mapsto v_c$ and which should be compared with \eqref{218equ} in the weight sequence setting:
$$v_c(t)\int_0^t\frac{1}{v_c(r)}dr=v(ct)\int_0^t\frac{1}{v(cr)}dr=v(ct)\int_0^{ct}\frac{1}{v(s)}\frac{ds}{c}.$$

The second supplement stated in the above result gives the next immediate consequence:

\begin{corollary}\label{Intopthmcor}
Let $\mathbf{M}$ be a log-convex weight sequence. If $\mathbf{M}$ is even strongly log-convex, then the integration operator $J: H^{\infty}_{v_{\mathbf{M}}}(\CC)\rightarrow H^{\infty}_{v_{\mathbf{M}}}(\CC)$ is not compact.
\end{corollary}

By combining Theorems \ref{Diffopthm} and \ref{Intopthm} we obtain:

\begin{corollary}\label{DiffIntopthmcor}
Let $\mathbf{M}$ be a log-convex weight sequence.
\begin{itemize}
	\item[$(i)$] The following are equivalent:
\begin{itemize}
\item[$(a)$] The differentiation operator $D: H^{\infty}_{v_{\mathbf{M}}}(\CC)\rightarrow H^{\infty}_{v_{\mathbf{M}}}(\CC)$ and the integration operator $J: H^{\infty}_{v_{\mathbf{M}}}(\CC)\rightarrow H^{\infty}_{v_{\mathbf{M}}}(\CC)$ are simultaneously continuous.

\item[$(b)$] $\mathbf{M}\hyperlink{approx}{\approx}\mathbf{G}^1$.
\end{itemize}
Moreover, when these equivalent conditions hold then neither $D$ nor $J$ is compact and in view of Lemma \ref{MG1equivlemma} the list of equivalent assertions can be extended by $(ii)-(v)$ in this auxiliary result.

\item[$(ii)$] The differentiation operator $D: H^{\infty}_{v_{\mathbf{M}}}(\CC)\rightarrow H^{\infty}_{v_{\mathbf{M}}}(\CC)$ and the integration operator $J: H^{\infty}_{v_{\mathbf{M}}}(\CC)\rightarrow H^{\infty}_{v_{\mathbf{M}}}(\CC)$ can never be simultaneously compact.
\end{itemize}
\end{corollary}

Note that assertion $(i)(b)$ is preserved under equivalence of sequences and $(ii)$ follows for abstractly weight functions directly by (b)(iii) in \cite[Thm. 26 \& 28]{Bonet2022survey}.

\demo{Proof}
$(i)(a)\Rightarrow(i)(b)$ This is an immediate consequence of the characterizations in Theorems \ref{Diffopthm} and \ref{Intopthm}; specifically see $(iv)\Rightarrow(v)$ in Theorem \ref{Intopthm}.

$(i)(b)\Rightarrow(i)(a)$ The equivalence implies $(iv)$ in Theorem \ref{Diffopthm} and $(v)$ in Theorem \ref{Intopthm}. Equivalence preserves \hyperlink{mg}{$(\on{mg})$}, hence this property follows for $\mathbf{M}$ since it is clear for $\mathbf{G}^1$. Thus $(v)\Rightarrow(iv)$ in Theorem \ref{Intopthm} is valid and we are done. Concerning the supplement note that $\mathbf{M}\hyperlink{approx}{\approx}\mathbf{G}^1$ obviously contradicts both $\mathbf{M}\hyperlink{mtriangle}{\vartriangleleft}\mathbf{G}^1$ and $\mathbf{G}^1\hyperlink{mtriangle}{\vartriangleleft}\mathbf{M}$ and then confirm $(iv)'$ in Theorem \ref{Diffopthm} and assumption $(v)'$ together with $(i)'\Rightarrow(v)'$  in Theorem \ref{Intopthm}.\vspace{6pt}

$(ii)$ This is immediate since $(v)'$ in Theorem \ref{Diffopthm} and $(iv)'$ in Theorem \ref{Intopthm} can never be satisfied simultaneously.
\qed\enddemo

Next, we investigate \cite[Cor. 29]{Bonet2022survey} (resp. \cite[Cor. 3.10]{AbaninThien17}) in the weight sequence setting.

\begin{theorem}\label{Intopthmcor}
	Let $\mathbf{M}$ be a log-convex weight sequence. Then the following are equivalent:
	\begin{itemize}
		\item[$(i)$] The differentiation operator $D: H^{\infty}_{v_{\mathbf{M}}}(\CC)\rightarrow H^{\infty}_{v_{\mathbf{M}}}(\CC)$ is continuous and surjective.
		
\item[$(ii$)] $0<\liminf_{r\rightarrow+\infty}v_{\mathbf{M}}(r)\left(\frac{1}{v_{\mathbf{M}}}\right)'(r)\le\limsup_{r\rightarrow+\infty}v_{\mathbf{M}}(r)\left(\frac{1}{v_{\mathbf{M}}}\right)'(r)<+\infty$.
		
\item[$(iii)$] $\mathbf{M}\hyperlink{approx}{\approx}\mathbf{G}^1$ holds.
	\end{itemize}

The list of equivalences can be extended by the assertions in Lemma \ref{MG1equivlemma} and $(iii)$ is preserved under equivalence of sequences.
\end{theorem}

\demo{Proof}
$(i)\Leftrightarrow(ii)$ holds by \cite[Cor. 29]{Bonet2022survey} applied to $v=v_{\mathbf{M}}$.

Concerning $(ii)\Leftrightarrow(iii)$ we give two arguments: First, by $(i)\Leftrightarrow(v)$ in Lemma \ref{MG1equivlemma} we immediately get $(iii)\Leftrightarrow(i)$ since $(v)$ in this lemma is precisely \cite[Cor. 29 (iii)]{Bonet2022survey}.

Second, also in view of Lemma \ref{MG1equivlemma}, in order to prove $(ii)\Leftrightarrow(iii)$ it suffices to show that $(ii)$ is equivalent to
$$0<\inf_{p\in\NN_{>0}}\frac{\mu_p}{p}\le\sup_{p\in\NN_{>0}}\frac{\mu_p}{p}<+\infty,$$
and by $(iii)\Leftrightarrow(iv)$ in Theorem \ref{Intopthm} for this it suffices to verify that $\limsup_{t\rightarrow+\infty}v_{\mathbf{M}}(t)\left(\frac{1}{v_{\mathbf{M}}}\right)'(t)<+\infty$ if and only if $0<\inf_{p\in\NN_{>0}}\frac{\mu_p}{p}$:\vspace{6pt}

If $\limsup_{t\rightarrow+\infty}v_{\mathbf{M}}(t)\left(\frac{1}{v_{\mathbf{M}}}\right)'(t)<+\infty$, then in view of the computation in the proof of Theorem \ref{Intopthm}:
$$\exists\;C\ge 1\;\forall\;p\in\NN_{>0}\;\forall\;t\ge 0,\;\mu_p\le t<\mu_{p+1}:\;\;\;\frac{p}{t}\le A.$$
Thus, the choice $t=\mu_p$ gives $\frac{\mu_p}{p}\ge A^{-1}$ for all $p\in\NN_{>0}$ such that $\mu_p<\mu_{p+1}$. And if $\mu_p=\dots=\mu_{p+\ell}<\mu_{p+\ell+1}$ for some $\ell\in\NN_{>0}$, then $\frac{\mu_{p+\ell}}{p+\ell}<\frac{\mu_{p+i}}{p+i}$ for all $0\le i\le\ell-1$ and so the first half is shown.

Conversely, $0<\inf_{p\in\NN_{>0}}\frac{\mu_p}{p}$ precisely means that $\epsilon<\frac{\mu_p}{p}$ for some $\epsilon>0$ and all $p\in\NN_{>0}$ which implies $\frac{p}{t}<\epsilon^{-1}$ for all $p\in\NN_{>0}$ and for all $t\in[\mu_p,\mu_{p+1})$ and so we are done.
\qed\enddemo

On the other hand, when involving the information from Lemma \ref{MuG1equivlemma} for which we recall that \cite[Cor. 29 (iii)]{Bonet2022survey} is $(v)$ in Lemma \ref{MuG1equivlemma} for the weight $u$, we can extend resp. reformulate \cite[Cor. 29]{Bonet2022survey}:

\begin{theorem}\label{Intopthmcoru}
Let $u$ be a normalized and (log-)convex weight function. Then the following are equivalent:
\begin{itemize}
\item[$(i)$] The differentiation operator $D: H^{\infty}_{u}(\CC)\rightarrow H^{\infty}_{u}(\CC)$ is continuous and surjective.
		
\item[$(ii$)] $0<\liminf_{t\rightarrow+\infty}u(t)\left(\frac{1}{u}\right)'(t)\le\limsup_{t\rightarrow+\infty}u(t)\left(\frac{1}{u}\right)'(t)<+\infty$.
		
\item[$(iii)$] It holds that
$$\exists\;C,C_1\ge 1\;\forall\;t\ge 0:\;\;\;\frac{1}{C_1}\exp(-Ct)\le u(t)\le C_1\exp(-t/C).$$

\item[$(iv)$] $\mathbf{M}^u\hyperlink{approx}{\approx}\mathbf{G}^1$ holds.
	\end{itemize}

The list of equivalences can be extended by the assertions listed in Lemma \ref{MuG1equivlemma}.
\end{theorem}

When combining Corollary \ref{Diffopthmcor}, $(i)$ in Corollary \ref{DiffIntopthmcor}, Theorems \ref{Intopthmcor} and \ref{Intopthmcoru}, and \cite[Thm. 26 \& 28, Cor. 29]{Bonet2022survey} the next result is shown:

\begin{corollary}\label{finalDiffIntcor}
We get:
\begin{itemize}
\item[$(I)$] Let $\mathbf{M}$ be a log-convex weight sequence. Then the following are equivalent:
\begin{itemize}
\item[$(i)$] The differentiation operator $D: H^{\infty}_{v_{\mathbf{M}}}(\CC)\rightarrow H^{\infty}_{v_{\mathbf{M}}}(\CC)$ is continuous and surjective.
		
\item[$(ii)$] The differentiation operator $D: H^{\infty}_{v_{\mathbf{M}}}(\CC)\rightarrow H^{\infty}_{v_{\mathbf{M}}}(\CC)$ and the integration operator $J: H^{\infty}_{v_{\mathbf{M}}}(\CC)\rightarrow H^{\infty}_{v_{\mathbf{M}}}(\CC)$ are simultaneously continuous.
\item[$(iii)$] $\mathbf{M}\hyperlink{approx}{\approx}\mathbf{G}^1$ holds.
	\end{itemize}
If any of these equivalent conditions hold then neither $D$ nor $J$ is compact.

\item[$(II)$] Let $u$ be a normalized and (log-)convex weight function. Then the following are equivalent:
\begin{itemize}
\item[$(i)$] The differentiation operator $D: H^{\infty}_{u}(\CC)\rightarrow H^{\infty}_{u}(\CC)$ is continuous and surjective.
		
\item[$(ii)$]  The differentiation operator $D: H^{\infty}_{u}(\CC)\rightarrow H^{\infty}_{u}(\CC)$ and the integration operator $J: H^{\infty}_{u}(\CC)\rightarrow H^{\infty}_{u}(\CC)$ are simultaneously continuous.

\item[$(iii)$] $\mathbf{M}^u\hyperlink{approx}{\approx}\mathbf{G}^1$ holds.

\item[$(iv)$] The differentiation operator $D: H^{\infty}_{v_{\mathbf{M}^u}}(\CC)\rightarrow H^{\infty}_{v_{\mathbf{M}^u}}(\CC)$ is continuous and surjective.
		
\item[$(v$)] The differentiation operator $D: H^{\infty}_{v_{\mathbf{M}^u}}(\CC)\rightarrow H^{\infty}_{v_{\mathbf{M}^u}}(\CC)$ and the integration operator $J: H^{\infty}_{v_{\mathbf{M}^u}}(\CC)\rightarrow H^{\infty}_{v_{\mathbf{M}^u}}(\CC)$ are simultaneously continuous.
	\end{itemize}
\end{itemize}
The equivalences listed in (I) can be extended by the assertions in Lemma \ref{MG1equivlemma}, whereas those in $(II)$ by the ones in Lemma \ref{MuG1equivlemma}.
\end{corollary}

In the above results $\mathbf{M}\hyperlink{approx}{\approx}\mathbf{G}^{1}$ appears crucially; for concrete examples this can be checked directly and we investigate this condition in more detail. First, in view of \eqref{equ14} and the comments in Section \ref{weightsequencesection} it is straight-forward to see that there exist log-convex weight sequences $\mathbf{M}$ such that $\mathbf{M}\hyperlink{approx}{\approx}\mathbf{G}^1$ holds but such that either
\begin{equation*}\label{counterrelationGevrey1}
\exists\;A\ge 1\;\forall\;p\in\NN:\;\;\;M_p\le AG^1_p=Ap!
\end{equation*}
or
\begin{equation*}\label{counterrelationGevrey2}
\exists\;A\ge 1\;\forall\;p\in\NN:\;\;\;p!=G^1_p\le AM_p
\end{equation*}
fails. Indeed, for this we can take $M_p:=h^pp!$, $0<h\neq 1$ fixed. Moreover, the characterization established in \cite[Thm. 3.5]{weightedentireinclusion1} gives that in this situation either $H^{\infty}_{v_{\mathbf{G}^1}}(\CC)\subseteq H^{\infty}_{v_{\mathbf{M}}}(\CC)$ or $H^{\infty}_{v_{\mathbf{M}}}(\CC)\subseteq H^{\infty}_{v_{\mathbf{G}^1}}(\CC)$ is violated and hence, in particular, the equality (as Banach spaces) $H^{\infty}_{v_{\mathbf{M}}}(\CC)=H^{\infty}_{v_{\mathbf{G}^1}}(\CC)$ fails. In this sense the equivalence relation \hyperlink{approx}{$\approx$} in the above results gives more flexibility; for example by $(i)$ in Corollary \ref{DiffIntopthmcor} the differentiation and integration operator are jointly continuous even for weighted spaces $H^{\infty}_{v_{\mathbf{M}}}(\CC)$ which do not coincide (as sets and/or Banach spaces) with $H^{\infty}_{v_{\mathbf{G}^1}}(\CC)$.

We strengthen these observations in the next technical statement which also illustrates the difference between the crucial characterizing growth relations in the weighted settings; recall \cite[Prop. 2.12 \& 4.6]{compositionpaper} in the ultradifferentiable weight sequence setting and \cite[Thm. 3.5]{weightedentireinclusion1}, \cite[Thm. 3.2]{weightedentireinclusion2}, Theorem \ref{Thm35disk} in the weighted analytic case.

\begin{lemma}\label{weakerrelationex}
Let $\mathbf{N}$ be a log-convex weight sequence. Then there exists a log-convex weight sequence $\mathbf{M}$ such that $\mathbf{M}\hyperlink{approx}{\approx}\mathbf{N}$ but such that both
\begin{equation}\label{counterrelation1}
\exists\;A\ge 1\;\forall\;p\in\NN:\;\;\;M_p\le AN_p
\end{equation}
and
\begin{equation}\label{counterrelation2}
\exists\;A\ge 1\;\forall\;p\in\NN:\;\;\;N_p\le AM_p
\end{equation}
are violated.
\end{lemma}

This result applied to $\mathbf{N}=\mathbf{G}^1$ and the above comments imply that for the constructed sequence $\mathbf{M}$ neither $H^{\infty}_{v_{\mathbf{M}}}(\CC)\subseteq H^{\infty}_{v_{\mathbf{G}^1}}(\CC)$ nor $H^{\infty}_{v_{\mathbf{G}^1}}(\CC)\subseteq H^{\infty}_{v_{\mathbf{M}}}(\CC)$ holds but both operators $D,J: H^{\infty}_{v_{\mathbf{M}}}(\CC)\rightarrow H^{\infty}_{v_{\mathbf{M}}}(\CC)$ are continuous and $D$ is surjective (see Theorem \ref{Intopthmcor}).

\demo{Proof}
For any $h>0$ write $\mathbf{N}^h:=(h^pN_p)_{p\in\NN}$ and generally note that $N^h_0=N_0$, $\mathbf{N}^h$ is a (log-convex) weight sequence (according to Definition \ref{defweightsequ}) if and only if $\mathbf{N}$ is so. Indeed, the normalization condition $M_0=1$ is not required in the forthcoming construction. By definition one has $\nu^h_p=h\nu_p$ for all $p\in\NN$ and $h>0$ (again with $\nu^h_0:=1$).

Fix now $h>1$ and $\mathbf{M}$ is defined in terms of the quotients $\mu=(\mu_p)_{p\in\NN_{>0}}$ via $M_p:=M_0\prod_{i=1}^p\mu_i$ and strictly increasing sequences of integers $\mathbf{a}:=(a_j)_{j\in\NN_{>0}}$, $\mathbf{b}:=(b_j)_{j\in\NN_{>0}}$, $\mathbf{c}:=(c_j)_{j\in\NN_{>0}}$ with $a_j<b_j\le c_j<a_{j+1}$ for all $j\in\NN_{>0}$ and such that $a_1:=1$. We put
\begin{equation}\label{weakerrelationexequ}
M_0:=N_0,
\end{equation}
\begin{equation}\label{weakerrelationexequ1}
\mu_p:=h\nu_p(=\nu^h_p),\hspace{15pt}a_j\le p\le b_j-1,\;j\in\NN_{>0},
\end{equation}
\begin{equation}\label{weakerrelationexequ2}
\mu_p:=h\nu_{b_j-1}(=\nu^h_{b_j-1}),\hspace{15pt}b_j\le p\le c_j-1,\;j\in\NN_{>0},
\end{equation}
\begin{equation}\label{weakerrelationexequ3}
\mu_p:=h^{-1}\nu_{p}(=\nu^{h^{-1}}_p),\hspace{15pt}c_j\le p\le a_{j+1}-1,\;j\in\NN_{>0}.
\end{equation}
First, by \eqref{weakerrelationexequ} it is immediate that $M_0=1$ if and only if $N_0=1$. Second, let $c_j$ be chosen minimal such that $h\nu_{b_j-1}<h^{-1}\nu_{c_j}\Leftrightarrow h^2\nu_{b_j-1}<\nu_{c_j}$, then by definition $p\mapsto\mu_p$ is non-decreasing and thus $\mathbf{M}$ is log-convex. And this choice for $c_j$ is possible since $\lim_{p\rightarrow+\infty}\nu_p=\lim_{p\rightarrow+\infty}(N_p)^{1/p}=+\infty$ holds.

Moreover, by log-convexity of $\mathbf{N}$ it follows that $c_j>b_j-1$ and hence \eqref{weakerrelationexequ1}, \eqref{weakerrelationexequ2} and \eqref{weakerrelationexequ3} and the minimality of $c_j$ ensure that $\nu^{h^{-1}}_p\le\mu_p\le\nu^h_p$ for all $p\in\NN_{>0}$. This gives $\mathbf{N}^{h^{-1}}\le\mathbf{M}\le\mathbf{N}^h$ and so $\mathbf{M}\hyperlink{approx}{\approx}\mathbf{N}$ is shown. Note that, when $c_j=b_j$ and which is not excluded in the general case, then the set of indices $p$ in \eqref{weakerrelationexequ2} is empty.

On the other hand, let us check that when (iteratively) $a_j$ and $b_j$ are chosen to be increasing sufficiently fast, then both \eqref{counterrelation1} and \eqref{counterrelation2} fail.

More precisely, by \eqref{weakerrelationexequ3} we get
$$M_{a_{j+1}-1}=M_{c_j-1}h^{-(a_{j+1}-c_j)}\frac{N_{a_{j+1}-1}}{N_{c_j-1}}\le\frac{1}{j}N_{a_{j+1}-1}\Leftrightarrow j\frac{M_{c_j-1}}{N_{c_j-1}}\le h^{a_{j+1}-c_j},$$
which can be achieved when choosing $a_{j+1}>c_j$ sufficiently large (recall that $h>1$ is valid). Consequently,
$$\forall\;j\in\NN_{>0}:\;\;\;jM_{a_{j+1}-1}\le N_{a_{j+1}-1},$$
and which obviously contradicts \eqref{counterrelation2}.

Finally, \eqref{weakerrelationexequ1} implies
$$M_{b_j-1}=M_{a_j-1}h^{b_j-a_j}\frac{N_{b_j-1}}{N_{a_j-1}}\ge jN_{b_j-1}\Leftrightarrow h^{b_j-a_j}\ge j\frac{N_{a_j-1}}{M_{a_j-1}},$$
which holds when $b_j>a_j$ is sufficiently large (again recall that $h>1$). This estimate implies
$$\forall\;j\in\NN_{>0}:\;\;\;M_{b_j-1}\ge jN_{b_j-1},$$
and thus \eqref{counterrelation1} fails, too.

Note that in the above arguments the sequence $(j)_{j\in\NN_{>0}}$ can be replaced by any arbitrary increasing sequence $\mathbf{d}:=(d_j)_{j\in\NN_{0}}$ satisfying $\lim_{j\rightarrow+\infty}d_j=+\infty$ when modifying $\mathbf{a}$, $\mathbf{b}$ and $\mathbf{c}$ accordingly.
\qed\enddemo

On the other hand, in the next step we show that log-convex weight sequences $\mathbf{M}$ satisfying $\mathbf{M}\hyperlink{approx}{\approx}\mathbf{G}^1$ can be constructed explicitly in a straight-forward and convenient way. Let $$u_{\alpha}(t):=\exp(-\alpha t),\;\;\;t\in[0,+\infty),\;\alpha>0,$$
which is standard in the theory of weighted entire spaces. The corresponding \emph{associated weight sequence} $\mathbf{M}^{\alpha}:=(M^{\alpha}_p)_{p\in\NN}$ given by
\begin{equation}\label{Malphadef}
M^{\alpha}_p:=\sup_{t\ge 0}t^pu_{\alpha}(t),\;\;\;p\in\NN,
\end{equation}
is equivalent to $\mathbf{G}^1$. More precisely, for $p=0$ one gets $M^{\alpha}_0=1$ and for $p\in\NN_{>0}$ put $f_p(t):=t^pu_{\alpha}(t)$, $t\in(0,+\infty)$. Then $f'_p(t)=pt^{p-1}\exp(-\alpha t)+t^p\exp(-\alpha t)(-\alpha)=0\Leftrightarrow p=t\alpha\Leftrightarrow\frac{p}{\alpha}=t$ and at this point the maximum of $f_p$ on $(0,+\infty)$ is attained and so (with $0^0:=1$):
\begin{equation}\label{Malpha}
\forall\;p\in\NN:\;\;\;M^{\alpha}_p=f_p\left(\frac{p}{\alpha}\right)=\left(\frac{p}{\alpha}\right)^pe^{-p}=p^p\left(\frac{1}{e\alpha}\right)^p.
\end{equation}
Thus $\mathbf{M}^{\alpha}\hyperlink{approx}{\approx}\mathbf{G}^{1}$ for any $\alpha>0$, more precisely by \eqref{stirling} (and $0^0:=1$):
\begin{equation}\label{Malphadef1}
\forall\;\alpha>0\;\forall\;p\in\NN:\;\;\;\left(\frac{1}{e\alpha}\right)^pG^1_p=p!(e\alpha)^{-p}\le p^p(e\alpha)^{-p}=M^{\alpha}_p\le p!\left(\frac{1}{\alpha}\right)^p=G^1_p\left(\frac{1}{\alpha}\right)^p.
\end{equation}
This computation should be compared with the proof of $(iii)\Rightarrow(iv)$ in Theorem \ref{Diffopthm} and note that $\mathbf{M}^{\alpha}$ is a log-convex weight sequence: Apart from $M^{\alpha}_0=1$, by definition it holds that $\mathbf{M}^{\alpha}$ is log-convex (cf. comment $(b)$ in Section \ref{assosequencesection}) and since $u_{\alpha}$ is strictly positive also $(M^{\alpha}_p)^{1/p}\rightarrow+\infty$ (recall \cite[Rem. 2.14]{weightedentireinclusion1} and $(d)$ in Section \ref{assosequencesection}). $u_{\alpha}$ is (log-)convex since $t\mapsto\varphi_{\omega^{u_{\alpha}}}(t)=\alpha e^t$ is obviously convex. $u_{\alpha}$ is not normalized but by Remark \ref{normrem} and Lemmas \ref{normconvexrem}, \ref{auxlemma} there exists a normalized and (log-)convex weight $u^n_{\alpha}$ such that $u^{n}_{\alpha}(t)=u_{\alpha}(t)$ for all $t$ sufficiently large which gives
\begin{equation}\label{ualphanvsualphaequ}
\forall\;\alpha>0\;\exists\;A_1\ge 1\;\forall\;t\ge 0:\;\;\;\frac{1}{A_1}u_{\alpha}(t)\le u^{n}_{\alpha}(t)\le A_1u_{\alpha}(t).
\end{equation}
Hence, by definitions \eqref{vBMTweight1equ1} and \eqref{Malphadef}, it follows that
\begin{equation}\label{Malphauequ}
\forall\;\alpha>0\;\exists\;A_1\ge 1\;\forall\;p\in\NN:\;\;\;\frac{1}{A_1}M^{\alpha}_p\le M^{u^n_{\alpha}}_p\le A_1M^{\alpha}_p.
\end{equation}
Then even $\mathbf{M}^{u^n_{\alpha}}\in\hyperlink{LCset}{\mathcal{LC}}$ and \eqref{Malphadef1}, \eqref{Malphauequ} give $\mathbf{M}^{u^n_{\alpha}}\hyperlink{approx}{\approx}\mathbf{G}^{1}$ for any $\alpha>0$, too.\vspace{6pt}

We close this section by commenting on a technical issue in \cite{orliczpaper}.

\begin{remark}\label{orliczpaperrem}
Motivated by the proof of $(iii)\Rightarrow(iv)$ in Theorem \ref{Intopthm} we see that the technical assumption $\mu_p<\mu_{p+1}$ for all $p\in\NN_{>0}$ in \cite[Lemma 4.9]{orliczpaper} is superfluous: Indeed, the given first part of $(i)\Rightarrow(ii)$ on \cite[p. 938]{orliczpaper} implies $\mu_j<\lambda_{\lceil\frac{j}{c_1}\rceil}$ for all $j\ge j_{\epsilon_1}$ such that $\mu_j<\mu_{j+1}$ and w.l.o.g. $\mu_{j_{\epsilon_1}}<\mu_{j_{\epsilon_1}+1}$. If now $j\ge j_{\epsilon_1}$ is such that $\mu_j<\mu_{j+1}=\dots=\mu_{j+1+\ell}<\mu_{j+2+\ell}$ for some $\ell\in\NN_{>0}$, then indeed the arguments yield $\Sigma_{\mathbf{L}}(t)<\lceil\frac{j}{c_1}\rceil$ and hence $t<\lambda_{\lceil\frac{j}{c_1}\rceil}$ for all $t\in[\mu_j,\mu_{j+1})$. Thus $\mu_{j+1}\le\lambda_{\lceil\frac{j}{c_1}\rceil}$ holds which implies by log-convexity of $\mathbf{L}$ that $\mu_{j+1+i}\le\lambda_{\lceil\frac{j+i}{c_1}\rceil}$ for all $i=0,\dots,\ell$. Summarizing, $\mu_j\le\lambda_{\lceil\frac{j}{c_1}\rceil}$ is valid for all $j\ge j_{\epsilon_1}$.

For the second part in this implication, when $\mu_j=\dots=\mu_{j+\ell}<\mu_{j+\ell+1}$, then one can choose $t=\mu_{j+\ell}$ and conclude again for all indices $j+i$, $i=0,\dots,\ell$, by the log-convexity of $\mathbf{L}$.

The rest of the proof of \cite[Lemma 4.9]{orliczpaper} does not require any changes.
\end{remark}

\subsection{The weighted entire case with $o$-growth restriction}\label{diffintsectionogrowth}
From now on we focus on weighted spaces of entire functions with $o$-growth restriction (at infinity). Again the sequence $\mathbf{G}^1$ plays a crucial role in our considerations. First, we transfer \cite[Thm. 33]{Bonet2022survey} resp. \cite[Thm. 2.1]{Bonet09} to the weight sequence setting.

\begin{theorem}\label{Thm33}
Let $\mathbf{M}$ be a log-convex weight sequence and assume that the differentiation operator $D: H^0_{v_{\mathbf{M}}}(\CC)\rightarrow H^0_{v_{\mathbf{M}}}(\CC)$ is continuous. Consider the following assertions:
\begin{itemize}
\item[$(i)$] $D$ has a dense set of periodic points.

\item[$(ii)$] $D$ has a periodic point different of $0$.

\item[$(iii)$] It holds that
\begin{equation}\label{Thm33equ}
\lim_{t\rightarrow+\infty}e^tv_{\mathbf{M}}(t)=0\Longleftrightarrow\lim_{r\rightarrow+\infty}\omega_{\mathbf{M}}(t)-t=+\infty.
\end{equation}
\item[$(iv)$] $t=o(\omega_{\mathbf{M}}(t))$ as $t\rightarrow+\infty$.

\item[$(v)$] $\mathbf{M}\hyperlink{mtriangle}{\vartriangleleft}\mathbf{G}^{1}$ is valid.

\item[$(vi)$] $\mathbf{M}\hyperlink{preceq}{\preceq}\mathbf{G}^{1}$ is valid.
\end{itemize}
Then we have the following implications/equivalences: $(i)\Leftrightarrow(ii)\Leftrightarrow(iii)\Rightarrow(vi)$ and $(iv)\Leftrightarrow(v)\Rightarrow(iii)$.
\end{theorem}

$(v)$ and $(vi)$ are obviously preserved under equivalence of (weight) sequences.

\demo{Proof}
$(i)\Leftrightarrow(ii)\Leftrightarrow(iii)$ holds by \cite[Thm. 34]{Bonet2022survey} (resp. \cite[Thm. 2.1]{Bonet09}) applied to $v=v_{\mathbf{M}}$.

And the equivalence stated in assertion $(iii)$ follows by $e^tv_{\mathbf{M}}(t)=\exp(t-\omega_{\mathbf{M}}(t))$; recall \eqref{weights}.\vspace{6pt}

$(iii)\Rightarrow(vi)$ \eqref{Thm33equ} implies, in particular, that
$$\exists\;D\in\RR\;\forall\;t\ge 0:\;\;\;\omega_{\mathbf{M}}(t)-t\ge D.$$
Following then the arguments from the proof of $(iii)\Rightarrow(iv)$ in Theorem \ref{Diffopthm} we see that this estimate gives
$$\exists\;D\in\RR\;\forall\;p\in\NN:\;\;\;M_p\le e^{-D}\frac{p^p}{e^p}\le e^{-D}p!,$$
which implies that $\mathbf{M}\hyperlink{preceq}{\preceq}\mathbf{G}^{1}$.\vspace{6pt}

$(iv)\Rightarrow(iii)$ If $t=o(\omega_{\mathbf{M}}(t))$ as $t\rightarrow+\infty$, i.e.
$$\forall\;\epsilon>0\;\exists\;C\ge 1\;\forall\;t\ge 0:\;\;\;t\le\epsilon\omega_{\mathbf{M}}(t)+C,$$
then $t-\omega_{\mathbf{M}}(t)\le(\epsilon-1)\omega_{\mathbf{M}}(t)+C$ for all $t\ge 0$. Choose a fixed $\epsilon\in(0,1)$ and then as $t\rightarrow+\infty$ one infers $\lim_{t\rightarrow+\infty}t-\omega_{\mathbf{M}}(t)=-\infty$; i.e. \eqref{Thm33equ}.\vspace{6pt}

Finally, $(iv)\Leftrightarrow(v)$ holds by repeating the arguments from $(iii)\Leftrightarrow(iv)$ resp. $(iii)'\Leftrightarrow(iv)'$ in Theorem \ref{Diffopthm}; see also \cite[Lemma 4.9 $(ii)\Leftrightarrow(iii)\Leftrightarrow(iv)$]{conjugateweightfunction}.
\qed\enddemo

\cite[Thm. 34]{Bonet2022survey} resp. \cite[Thm. 2.3]{Bonet09} takes the following form in the weight sequence setting:

\begin{theorem}\label{hypercyclthm34}
Let $\mathbf{M}$ be a log-convex weight sequence and assume that the differentiation operator $D: H^0_{v_{\mathbf{M}}}(\CC)\rightarrow H^0_{v_{\mathbf{M}}}(\CC)$ is continuous. Then the following are equivalent:
	\begin{itemize}
		\item[$(i)$] $D$ is \emph{hypercyclic} on $H^0_{v_{\mathbf{M}}}(\CC)$.
		
		\item[$(ii$)] It holds that $\liminf_{p\rightarrow+\infty}\frac{\|z^p\|_{v_{\mathbf{M}}}}{p!}=0$.
		
		\item[$(iii)$] $\mathbf{M}$ satisfies
\begin{equation}\label{hypercyclthmequ}
\liminf_{p\rightarrow+\infty}\frac{M_p}{p!}=\liminf_{p\rightarrow+\infty}m_p=0.
\end{equation}
\end{itemize}
\end{theorem}

\demo{Proof}
$(i)\Leftrightarrow(ii)$ holds by \cite[Thm. 34]{Bonet2022survey} applied to $v=v_{\mathbf{M}}$.

For $(ii)\Leftrightarrow(iii)$ note that by definition and \eqref{Prop32Komatsu}
$$\forall\;p\in\NN:\;\;\;\|z^p\|_{v_{\mathbf{M}}}=\sup_{z\in\CC}|z|^pv_{\mathbf{M}}(z)=\sup_{t\ge 0}\frac{t^p}{\exp(\omega_{\mathbf{M}}(t))}=M_p;$$
recall also \cite[$(2.26)$, Lemma 2.16]{weightedentireinclusion1}. Note that for the last equality log-convexity of $\mathbf{M}$ is indispensable.
\qed\enddemo

On the other hand, by recognizing that \cite[Thm. 2.3]{Bonet09} does not require necessarily (log-)convexity for $u$, using \eqref{vBMTweight1equ2} and applying the above proof to $\mathbf{M}^u$, we can immediately extend \cite[Thm. 34]{Bonet2022survey} resp. \cite[Thm. 2.3]{Bonet09} (a weight $v$ in \cite{Bonet09} satisfies the same conditions as in this work):

\begin{theorem}\label{hypercyclthm34fctversion}
Let $u$ be a weight function. Assume that the operator $D: H^0_{u}(\CC)\rightarrow H^0_{u}(\CC)$ is continuous. Then the following are equivalent:
	\begin{itemize}
		\item[$(i)$] $D$ is \emph{hypercyclic} on $H^0_{u}(\CC)$.
		
		\item[$(ii$)] It holds that $\liminf_{p\rightarrow+\infty}\frac{\|z^p\|_u}{p!}=0$.
		
		\item[$(iii)$] $\mathbf{M}^u$ satisfies \eqref{hypercyclthmequ}.
\end{itemize}
\end{theorem}
Note that (log-)convexity is not required for this result (see $(a)$ in Section \ref{assosequencesection}).\vspace{6pt}

Consider again $\mathbf{G}^s$, then \eqref{hypercyclthmequ} fails whenever $s\ge 1$ whereas for all $0<s<1$ this condition holds and also the (counter-)example in \cite[Prop. 6.9]{conjugateweightfunction} satisfies \eqref{hypercyclthmequ}. Indeed, any (log-convex weight) sequence $\mathbf{M}$ satisfying $\mathbf{M}\hyperlink{mtriangle}{\vartriangleleft}\mathbf{G}^{1}$ implies \eqref{hypercyclthmequ} but which is also non-standard in the ultradifferentiable setting. Moreover, we note:

\begin{remark}\label{hypercyclthmequrem}
\eqref{hypercyclthmequ} is preserved, when $\mathbf{M}$ and $\mathbf{N}$ are related by
$$\exists\;C\ge 1\;\forall\;p\in\NN:\;\;\;C^{-1}M_p\le N_p\le CN_p;$$
cf. \eqref{equ14}. However, it is not clear if \eqref{hypercyclthmequ} is preserved under equivalence of sequences and when $\mathbf{M}$ satisfies \eqref{hypercyclthmequ} then in general it is also not clear if some/any $\mathbf{M}^c$, $c\neq 1$, does have \eqref{hypercyclthmequ}; recall Section \ref{parametersection} for this notation and see Corollary \ref{hypercyclthmequlemmacor} and Examples \ref{example1} and \ref{example2}, Theorem \ref{hypercyclcor36} and Remark \ref{hypercyclthmequrem1}.
\end{remark}

We study now how \eqref{hypercyclthmequ} is related to the growth of $p\mapsto\frac{\mu_p}{p}$.

\begin{lemma}\label{hypercyclthmequlemma}
Let $\mathbf{M}\in\RR_{>0}^{\NN}$ be arbitrary. Then $\liminf_{p\rightarrow+\infty}\frac{\mu_p}{p}>1$ implies $\liminf_{p\rightarrow+\infty}\frac{M_p}{p!}>0$ and so \eqref{hypercyclthmequ} fails.
\end{lemma}

\demo{Proof}
By assumption there exists $p_1\in\NN_{>0}$ such that $\mu_p\ge p$ for all $p>p_1$ and which gives
$$\forall\;q>p_1:\;\;\;\frac{M_q}{q!}=\frac{M_0}{0!}\frac{\mu_1\cdots\mu_{p_1}}{1\cdots p_1}\frac{\mu_{p_1+1}\cdots\mu_q}{(p_1+1)\cdots q}\ge\frac{M_{p_1}}{p_1!}>0.$$
\qed\enddemo

\begin{corollary}\label{hypercyclthmequlemmacor}
Let $\mathbf{M}$ be a strongly log-convex weight sequence. If $\mathbf{M}$ satisfies \eqref{hypercyclthmequ}, then $\mathbf{M}\hyperlink{approx}{\approx}\mathbf{G}^1$.
\end{corollary}

\demo{Proof}
First, strong log-convexity gives that $\lim_{p\rightarrow+\infty}\frac{\mu_p}{p}\in(0,+\infty]$ since $p\mapsto\frac{\mu_p}{p}$ is non-decreasing.

However, by Lemma \ref{hypercyclthmequlemma} the value $+\infty$ for this limit is excluded since in this case $\lim_{p\rightarrow+\infty}\left(\frac{M_p}{p!}\right)^{1/p}=+\infty$ and which contradicts assumption \eqref{hypercyclthmequ}. Consequently, we have
$$\exists\;0<c<1<C<+\infty\;\forall\;p\in\NN_{>0}:\;\;\;c\le\frac{\mu_p}{p}\le C.$$
(Indeed, one can set $C:=\sup_{p\in\NN_{>0}}\frac{\mu_p}{p}$ and $c:=\mu_1$.) In view of Lemma \ref{MG1equivlemma} this is precisely $\mathbf{M}\hyperlink{approx}{\approx}\mathbf{G}^1$.
\qed\enddemo

Since $\mathbf{G}^1$ violates \eqref{hypercyclthmequ} it follows that in Corollary \ref{hypercyclthmequlemmacor} one cannot expect an equivalence between \eqref{hypercyclthmequ} and $\mathbf{M}\hyperlink{approx}{\approx}\mathbf{G}^1$; i.e. the latter condition is only necessary. This observation motivates the next (counter-)example which illustrates, together with Corollary \ref{hypercyclthmequlemmacor}, that \eqref{hypercyclthmequ} is in general not preserved under equivalence of (weight) sequences.

\begin{example}\label{example1}
There exists a strongly log-convex weight sequence $\mathbf{M}$ such that \eqref{hypercyclthmequ} holds.

Let $\mathbf{b}=(b_j)_{j\in\NN_{>0}}$ be a (fixed) strictly increasing sequence such that $b_j\in(0,1)$ for all $j\in\NN_{>0}$ and let $\mathbf{a}:=(a_j)_{j\in\NN_{>0}}$ be a strictly increasing sequence of integers with $a_1:=1$. Then define $\mathbf{M}$ by
\begin{equation}\label{Mdef}
\forall\;p\in\NN_{>0}:\;\;\;M_p:=M_0\prod_{i=1}^p\mu_i,\hspace{15pt}M_0:=1,
\end{equation}
with
\begin{equation}\label{example1mudef}
\mu_p:=b_jp,\;\;\;a_j\le p\le a_{j+1}-1,\;j\in\NN_{>0}.
\end{equation}
Obviously, $p\mapsto\frac{\mu_p}{p}$ is non-decreasing and $(0<)b_1\le\frac{\mu_p}{p}\le\sup_{j\in\NN_{>0}}b_j\le 1$ which gives strong log-convexity and the equivalence $\mathbf{M}\hyperlink{approx}{\approx}\mathbf{G}^1$ as a necessary condition for \eqref{hypercyclthmequ} via Corollary \ref{hypercyclthmequlemmacor}. Moreover, by this equivalence it is immediate that $\mathbf{M}$ is also a log-convex weight sequence according to Definition \ref{defweightsequ}. Finally, for all $j\in\NN_{>0}$ we obtain
$$\frac{M_{a_{j+1}-1}}{(a_{j+1}-1)!}=\prod_{i=1}^{a_{j+1}-1}\frac{\mu_i}{i}=\prod_{i=1}^jb_i^{a_{i+1}-a_i}\le\prod_{i=1}^jb_j^{a_{i+1}-a_i}=b_j^{a_{j+1}-a_1}=b_j^{a_{j+1}-1},$$
and when $\mathbf{a}$ is increasing sufficiently fast then, since $b_j\in(0,1)$ for all $j$, we can achieve that $b_j^{a_{j+1}-1}\le\frac{1}{j}$ for all $j\in\NN_{>0}$. This gives \eqref{hypercyclthmequ} and note that, when modifying $\mathbf{a}$ accordingly, $\frac{1}{j}$ can be replaced by $\frac{1}{c_j}$, with $\mathbf{c}=(c_j)_{j\in\NN_{>0}}$ denoting an arbitrary strictly increasing sequence tending to infinity; i.e. the rate of approaching $0$ (in \eqref{hypercyclthmequ}) is not effecting the validity of this statement.
\end{example}

A modification of this first example implies that even the modification $\mathbf{M}\mapsto\mathbf{M}^c$ is in general not preserving \eqref{hypercyclthmequ}.

\begin{example}\label{example2}
There exists a strongly log-convex weight sequence $\mathbf{M}$ such that \eqref{hypercyclthmequ} holds but for a precise choice $c>1$ the sequence $\mathbf{M}^c$, which is also strongly log-convex, violates \eqref{hypercyclthmequ}.\vspace{6pt}

We consider the sequence $\mathbf{M}$ defined in Example \ref{example1}, see \eqref{example1mudef}, and thus \eqref{hypercyclthmequ} holds. However, choose now $c:=\frac{1}{b_1}\in(1,+\infty)$ and consider the sequence $\mathbf{M}^c$. Hence, $\mathbf{M}^c$ is also a strong log-convex weight sequence being equivalent to $\mathbf{G}^1$, and in view of \eqref{example1mudef} we get for the corresponding sequence of quotients:
\begin{equation*}\label{example2mudef}
\mu^c_p:=cb_jp,\;\;\;a_j\le p\le a_{j+1}-1,\;j\in\NN_{>0}.
\end{equation*}
Therefore, since $b_j\ge b_1(>0)$ for all $j\in\NN_{>0}$ it holds that $\mu^c_p\ge p$ for all $p\in\NN_{>0}$ which implies $M^c_p\ge p!$ for all $p\in\NN$ (the case $p=0$ is clear by $M^c_0=M_0=1$). But this estimate obviously contradicts \eqref{hypercyclthmequ}.
\end{example}

On the other hand, we have the following (counter-)example which shows that \eqref{hypercyclthmequ} can even hold when $p\mapsto\frac{\mu_p}{p}$ is unbounded. Since by Lemma \ref{hypercyclthmequlemma} then necessarily $\liminf_{p\rightarrow+\infty}\frac{\mu_p}{p}\le 1$ also has to be satisfied, in this case the function $p\mapsto\frac{\mu_p}{p}$ has to be ``strongly oscillating''.

\begin{example}
There exists a log-convex weight sequence $\mathbf{M}=(M_p)_{p\in\NN}$ with \eqref{hypercyclthmequ} but such that
\begin{equation}\label{hypercyclthmequviol}
\sup_{p\in\NN_{>0}}\frac{\mu_p}{p}=+\infty.
\end{equation}
Let $\mathbf{a}=(a_j)_{j\in\NN_{>0}}$ be a strictly increasing sequence of integers with $a_1=1$ and let $\mathbf{b}=(b_j)_{j\in\NN_{>0}}$ be an arbitrary strictly increasing sequence in $\RR_{>0}$ such that $\lim_{j\rightarrow+\infty}b_j=+\infty$. Then introduce $\mathbf{M}$ again like in \eqref{Mdef} but with
$$\mu_p:=a_jb_j,\;\;\;\;\;\;a_j\le p\le a_{j+1}-1,\;j\in\NN_{>0}.$$
It is immediate that $p\mapsto\mu_p$ is non-decreasing and $\lim_{p\rightarrow+\infty}\mu_p=+\infty$, hence $\mathbf{M}$ is a weight sequence according to Definition \ref{defweightsequ} and log-convex. Moreover,
$$\sup_{p\in\NN_{>0}}\frac{\mu_p}{p}\ge \sup_{j\in\NN_{>0}}\frac{\mu_{a_j}}{a_j}=\sup_{j\in\NN_{>0}}b_j=+\infty$$
and so \eqref{hypercyclthmequviol} follows. However, if $\mathbf{a}$ is assumed to be growing fast enough, then also \eqref{hypercyclthmequ} holds: For this fix $\epsilon\in(0,1)$ and the auxiliary sequence $\mathbf{b}$, then when $a_j$ is given (start with $a_1:=1$) we can find $a'_j>a_j$ with $a'_j\le a_{j+1}-1$ such that $\frac{\mu_p}{p}=\frac{a_jb_j}{p}\le\epsilon$ for all $p>a'_j$. Thus, when choosing $a_{j+1}>a_j$ large enough one can ensure
\begin{equation}\label{hypercycliccounterex2equ}
\forall\;j\in\NN_{>0}:\;\frac{\mu_{a_j}\cdots\mu_{a_{j+1}-1}}{a_j\cdots(a_{j+1}-1)}=\frac{\mu_{a_j}\cdots\mu_{a'_j}}{a_j\cdots a'_j}\frac{\mu_{a'_j+1}\cdots\mu_{a_{j+1}-1}}{(a'_j+1)\cdots(a_{j+1}-1)}\le\frac{(a_jb_j)^{a'_j-a_j+1}}{a_j\cdots a'_j}\epsilon^{a_{j+1}-a'_j-1}\le\frac{1}{2},
\end{equation}
which implies, by setting $C_j:=\frac{\mu_{a_j}\cdots\mu_{a_{j+1}-1}}{a_j\cdots(a_{j+1}-1)}=\frac{(a_jb_j)^{a_{j+1}-a_j}}{a_j\cdots(a_{j+1}-1)}$, that $\frac{M_{a_{k+1}-1}}{(a_{k+1}-1)!}=\prod_{j=1}^kC_j\le\frac{1}{2^k}$ for all $k\in\NN_{>0}$. This proves \eqref{hypercyclthmequ} and concludes the construction. Note: Similarly as before, $\frac{1}{2}$ in \eqref{hypercycliccounterex2equ} can be replaced by any arbitrary strictly decreasing sequence $\mathbf{c}=(c_j)_{j\in\NN_{>0}}$ such that $\lim_{j\rightarrow+\infty}c_j=0$ when modifying $\mathbf{a}$ accordingly. Also the rate of unboundedness expressed in terms of $\mathbf{b}$ can be chosen as fast as desired.
\end{example}

\cite[Prop. 35]{Bonet2022survey} transfers immediately to the weight sequence setting; i.e. let $\mathbf{M}$ be a log-convex weight sequence, then $J$ is never hypercyclic on $H^0_{v_{\mathbf{M}}}(\CC)$ and has no periodic points different from $0$ in $H^{\infty}_{v_{\mathbf{M}}}(\CC)$.\vspace{6pt}

Consider now the function $u(t):=t^{1/2}e^{-t}$ on $[0,+\infty)$ crucially appearing in \cite[Cor. 36 (1)]{Bonet2022survey} and set accordingly
\begin{equation}\label{Nauxdef}
N_p:=\sup_{t\ge 0}t^pu(t)=\sup_{t\ge 0}t^{p+1/2}e^{-t},\;\;\;p\in\NN.
\end{equation}
Note that formally $u$ is not a weight function since, first, $u$ is not decreasing on whole $[0,+\infty)$ and, second, $u(0)=0$ and hence strict positivity on $[0,+\infty)$ also fails. Similarly as before let us compute $N_p$ explicitly: For $p\in\NN$ put $f_p(t):=t^{p+1/2}e^{-t}$, $t\in(0,+\infty)$, and so
$$f'_p(t)=(p+1/2)t^{p-1/2}e^{-t}+t^{p+1/2}(-1)e^{-t}=0\Leftrightarrow(p+1/2)=t,$$
and at this point the maximum of $f_p$ on $(0,+\infty)$ is attained. Thus for $\mathbf{N}=(N_p)_{p\in\NN}$ defined in \eqref{Nauxdef} we get
\begin{equation}\label{Nauxdef1}
\forall\;p\in\NN:\;\;\;N_p=f_p((p+1/2))=(p+1/2)^{p+1/2}e^{-p-1/2}.
\end{equation}
Indeed, since for the critical points $t_p:=p+1/2$ one has $t_p\ge 1/2$ for all $p\in\NN$ we see that $\mathbf{N}$ can be defined when considering the supremum in \eqref{Nauxdef} for all $t\ge t_0$ and $t_0\in(0,1/2)$ arbitrary but from now on fixed. Since $u$ vanishes at $0$ it (only) follows that
\begin{equation}\label{unvsuequ}
\forall\;t_1>0\;\exists\;A_2\ge 1\;\forall\;t\ge t_1:\;\;\;\frac{1}{A_2}u(t)\le u^{n}(t)\le A_2u(t),
\end{equation}
with $u^n$ denoting a normalized weight function (within the weighted entire setting) such that $u^n(t)=u(t)$ for all $t$ sufficiently large. The explicit construction of $u^n$ can be achieved similarly like in Remark \ref{normrem} (note that $u$ is rapidly decreasing). Moreover, since $u^n$ is normalized we can restrict in \eqref{vBMTweight1equ1} in the corresponding supremum defining $\mathbf{M}^{u^n}$ to all $t\ge t_2$ with $t_2\in[0,1]$ arbitrary. Summarizing, we apply \eqref{vBMTweight1equ1}, \eqref{Nauxdef}  and \eqref{unvsuequ} to $t_2=t_1=t_0\in(0,1/2)$ (fixed). These combined estimates yield analogously as in \eqref{assocsequnormrelation} and with the same $A_2$ as in \eqref{unvsuequ}:
\begin{equation}\label{MNuequ}
\exists\;A_2\ge 1\;\forall\;p\in\NN:\;\;\;\frac{1}{A_2}N_p\le M^{u^n}_p\le A_2N_p.
\end{equation}
By definition and, similarly as for $\mathbf{M}^{\alpha}$ above, we get that $\mathbf{N}$ is log-convex and $\mathbf{M}^{u^n}\in\hyperlink{LCset}{\mathcal{LC}}$.

Using this preparation and notation, \cite[Cor. 36]{Bonet2022survey} takes the following form:

\begin{theorem}\label{hypercyclcor36}
Let $\mathbf{M}$ be a log-convex weight sequence and assume that the differentiation operator $D: H^0_{v_{\mathbf{M}}}(\CC)\rightarrow H^0_{v_{\mathbf{M}}}(\CC)$ is continuous.
	\begin{itemize}
		\item[$(i)$] Consider the following assertions:
\begin{itemize}
\item[$(a)$] \begin{equation}\label{hypercyclcor36equ}
\exists\;A>0\;\forall\;t>0:\;\;\;\frac{1}{v_{\mathbf{M}}(t)}=\exp(\omega_{\mathbf{M}}(t))\le At^{-1/2}e^t=A\frac{1}{u(t)},
\end{equation}
\item[$(b)$]
\begin{equation}\label{hypercyclcor36equ0}
\forall\;t_0>0\;\exists\;A>0\;\forall\;t\ge t_0:\;\;\;\frac{1}{v_{\mathbf{M}}(t)}=\exp(\omega_{\mathbf{M}}(t))\le At^{-1/2}e^t=A\frac{1}{u(t)},
\end{equation}
\item[$(c)$] \begin{equation}\label{hypercyclcor36equ2}
\forall\;t_0>0\;\exists\;A>0\;\forall\;t\ge t_0:\;\;\;\frac{1}{v_{\mathbf{M}}(t)}=\exp(\omega_{\mathbf{M}}(t))\le A\frac{1}{u^n(t)},
\end{equation}
\item[$(d)$] \begin{equation}\label{hypercyclcor36equ3}
\exists\;A>0\;\forall\;p\in\NN:\;\;\;(p+1/2)^{p+1/2}\frac{1}{e^{p+1/2}}=N_p\le AM_p,
\end{equation}
\item[$(e)$]\begin{equation}\label{hypercyclcor36equ4}
\exists\;A>0\;\forall\;p\in\NN:\;\;\;M^{u^n}_p\le AM_p,
\end{equation}
\item[$(f)$] The operator $D$ is not hypercyclic on $H^0_{v_{\mathbf{M}}}(\CC)$.
\end{itemize}
Then $(a)\Leftrightarrow(b)\Leftrightarrow(c)\Leftrightarrow(d)\Leftrightarrow(e)\Rightarrow(f)$.

\item[$(ii)$] Consider the following assertions (and recall the definition of $\mathbf{M}^{\alpha}$ in \eqref{Malphadef}, \eqref{Malpha}):
\begin{itemize}
\item[$(a)$] \begin{equation}\label{hypercyclcor36equ5}
\exists\;B>0\;\exists\;\alpha>0\;\forall\;t\ge 0:\;\;\;v_{\mathbf{M}}(t)\le Bu_{\alpha}(t)=B\exp(-\alpha t),
\end{equation}
\item[$(b)$] \begin{equation}\label{hypercyclcor36equ6}
\exists\;B>0\;\exists\;\alpha>0\;\forall\;t\ge 0:\;\;\;v_{\mathbf{M}}(t)\le Bu^n_{\alpha}(t),
\end{equation}
\item[$(c)$] \begin{equation}\label{hypercyclcor36equ7}
\exists\;B>0\;\exists\;\alpha>0\;\forall\;p\in\NN:\;\;\;M_p\le BM^{\alpha}_p=Bp^p(e\alpha)^{-p},
\end{equation}
\item[$(d)$]\begin{equation}\label{hypercyclcor36equ8}
\exists\;B>0\;\exists\;\alpha>0\;\forall\;p\in\NN:\;\;\;M_p\le BM^{u_{\alpha}^n}_p,
\end{equation}
\item[$(e)$] The operator $D$ is hypercyclic on $H^0_{v_{\mathbf{M}}}(\CC)$.
\end{itemize}
Then $(a)\Leftrightarrow(b)\Leftrightarrow(c)\Leftrightarrow(d)$; if any of these assertions holds with $\alpha\ge 1$, then $(e)$ follows.
\end{itemize}
\end{theorem}

\emph{Note:} Apriori, the formulation in terms of \eqref{hypercyclcor36equ0} is natural since $\lim_{t\rightarrow 0}\frac{1}{u(t)}=+\infty$.

\demo{Proof}
$(i)$ \eqref{hypercyclcor36equ0}$\Leftrightarrow$\eqref{hypercyclcor36equ2} holds by \eqref{unvsuequ} and for \eqref{hypercyclcor36equ3}$\Leftrightarrow$\eqref{hypercyclcor36equ4} we involve \eqref{MNuequ}.

Choose in \eqref{hypercyclcor36equ2} a value $t_0\in(0,1]$ and then \eqref{Prop32Komatsu} and \eqref{vBMTweight1equ1} applied to $u^n$ immediately give \eqref{hypercyclcor36equ4}. On the other hand, \eqref{hypercyclcor36equ4} implies by definition $\omega_{\mathbf{M}}(t)\le\omega_{\mathbf{M}^{u^n}}(t)+\log(A)\Leftrightarrow v_{\mathbf{M}^{u^n}}(t)\le Av_{\mathbf{M}}(t)$ for all $t\ge 0$. The second half of \eqref{omegavequiv} applied to $u^n$, and for which (log-)convexity is not required, implies
$$\exists\;A>0\;\forall\;t\ge 0:\;\;\;u^n(t)\le v_{\mathbf{M}^{u^n}}(t)\le Av_{\mathbf{M}}(t),$$
and from this \eqref{hypercyclcor36equ2} follows immediately.

It remains to verify \eqref{hypercyclcor36equ}$\Leftrightarrow$\eqref{hypercyclcor36equ0} and the direction $\Rightarrow$ is clear. For the converse note that, indeed, $t\mapsto t^{-1/2}e^t$ is tending to infinity as $t\rightarrow 0^{+}$. Similarly as in the proof of \eqref{Nauxdef1}, one can show that $t\mapsto t^{-1/2}e^t$ is non-increasing on $(0,1/2]$. Also note that by \eqref{assovanishing} the left-hand side in \eqref{hypercyclcor36equ0} yields the constant value $1$ for all $t\in[0,\mu_1]$ and so, when choosing $t_0\in(0,\min\{\mu_1,1/2\})$ in \eqref{hypercyclcor36equ0}, then with the same choice for $A$ the estimate is valid for any $t>0$ and hence \eqref{hypercyclcor36equ}.

Thus all desired equivalences are shown and \eqref{hypercyclcor36equ} implies $(f)$ by \cite[Cor. 36 $(1)$]{Bonet2022survey} applied to $v=v_{\mathbf{M}}$.\vspace{6pt}

$(ii)$ \eqref{hypercyclcor36equ5}$\Leftrightarrow$\eqref{hypercyclcor36equ6} is valid by \eqref{ualphanvsualphaequ} and \eqref{hypercyclcor36equ7}$\Leftrightarrow$\eqref{hypercyclcor36equ8} holds by \eqref{Malphauequ}.

\eqref{hypercyclcor36equ6} implies \eqref{hypercyclcor36equ8} by \eqref{Prop32Komatsu} and \eqref{vBMTweight1equ1}. In all implications so far the value of the parameter $\alpha>0$ is preserved. Finally, \eqref{hypercyclcor36equ8} yields by definition $\omega_{\mathbf{M}^{u_{\alpha}^n}}(t)\le\omega_{\mathbf{M}}(t)+\log(B)\Leftrightarrow v_{\mathbf{M}}(t)\le Bv_{\mathbf{M}^{u_{\alpha}^n}}(t)$ for some $B>0$ and all $t\ge 0$. Then the first half of \eqref{omegavequiv} (applied to $u_{\alpha}^n$) and \eqref{ualphanvsualphaequ} give $v_{\mathbf{M}}(t)\le Bv_{\mathbf{M}^{u_{\alpha}^n}}(t)\le B\sqrt{A}\sqrt{u_{\alpha}^n(t)}\le B\sqrt{AA_1}\sqrt{u_{\alpha}(t)}=B\sqrt{AA_1}u_{\alpha/2}(t)$ for all $t\ge 0$ and hence \eqref{hypercyclcor36equ5} is shown with $\alpha'=\alpha/2$.

By these arguments we get that when we can choose $\alpha\ge 1$ in \eqref{hypercyclcor36equ5}, then $(e)$ follows by \cite[Cor. 36 $(2)$]{Bonet2022survey} applied to $v=v_{\mathbf{M}}$ and, moreover, if $\alpha\ge 2$ in \eqref{hypercyclcor36equ7} is possible then we can apply (again) \cite[Cor. 36 $(2)$]{Bonet2022survey} directly. However, note that \eqref{hypercyclcor36equ7} and \eqref{stirling} imply
$$\exists\;B\ge 1\;\forall\;p\in\NN:\;\;\;\frac{M_p}{p!}\le B\frac{p^p}{p!}(e\alpha)^{-p}\le B\alpha^{-p},$$
hence \eqref{hypercyclthmequ} holds for $\mathbf{M}^{\alpha}$ for any $\alpha>1$. Finally, if $\alpha=1$ then $\mathbf{M}^{\alpha}=(p^pe^{-p})_{p\in\NN}$ and this sequence satisfies \eqref{hypercyclthmequ} since by \emph{Stirling's formula} it holds that $p\mapsto\frac{p^p}{p!e^p}$ is asymptotically decreasing like $p\mapsto\frac{1}{\sqrt{2\pi p}}$ and which tends to $0$ as $p\rightarrow+\infty$.

Summarizing, $\mathbf{M}^{\alpha}$ satisfies \eqref{hypercyclthmequ} for any $\alpha\ge 1$ and hence, if \eqref{hypercyclcor36equ7} is valid for some $\alpha\ge 1$, then \eqref{hypercyclthmequ} for $\mathbf{M}$ follows and Theorem \ref{hypercyclthm34} yields the conclusion.
\qed\enddemo

Based on the arguments given in the second part of the above proof we point out (see also Remark \ref{hypercyclthmequrem}):

\begin{remark}\label{hypercyclthmequrem1}
First, also in $(i)$ the necessity of $(f)$ can be shown by involving Theorem \ref{hypercyclthm34}. For this note that the sequence $\mathbf{N}$ defined via \eqref{Nauxdef} violates \eqref{hypercyclthmequ} since for all $p\in\NN_{>0}$
$$\frac{N_p}{p!}=\frac{(p+1/2)^{1/2}}{e^{1/2}}\frac{(p+1/2)^p}{e^pp!}\ge\frac{p^{1/2}}{e^{1/2}}\frac{p^p}{e^pp!},$$
and again by \emph{Stirling's formula} it holds that $p\mapsto\frac{p^{1/2}}{e^{1/2}}\frac{p^p}{e^pp!}$ is approaching $\frac{1}{\sqrt{2\pi e}}$. Now, when $\mathbf{M}$ is given and satisfying \eqref{hypercyclcor36equ3}, then also \eqref{hypercyclthmequ} fails for $\mathbf{M}$.\vspace{6pt}

Second, concerning $(ii)$ in the above theorem in view of \eqref{stirling} assertion \eqref{hypercyclcor36equ7} precisely means that $\mathbf{M}\hyperlink{preceq}{\preceq}\mathbf{G}^{1}$ has to be satisfied and hence \eqref{hypercyclcor36equ7} is obviously preserved under equivalence of weight sequences. On the other hand, for the implication $(d)\Rightarrow(e)$ the parameter $\alpha$ plays the crucial role here but relation \hyperlink{approx}{$\approx$} \emph{does not recognize this subtle behavior.} More precisely, fix $\alpha_0>0$ in \eqref{hypercyclcor36equ7}, then each sequence $\mathbf{M}^{\alpha}$ (given by $M^{\alpha}_p:=p^p(e\alpha)^{-p}$) with $\alpha\ge\alpha_0$ satisfies \eqref{hypercyclcor36equ7} whereas this condition fails for all $\mathbf{M}^{\alpha}$ with $\alpha\in(0,\alpha_0)$. On the other hand, recall that obviously $\mathbf{M}^{\alpha}\hyperlink{approx}{\approx}\mathbf{M}^{\beta}(\hyperlink{approx}{\approx}\mathbf{G}^1)$ holds for all $\alpha,\beta>0$.
\end{remark}

Finally, in \cite[Thm. 39]{Bonet2022survey} the following (sufficient) condition on a weight function appears:

\begin{equation}\label{Thm39liminf}
\liminf_{t\rightarrow+\infty}\frac{tv_{\mathbf{M}}(t)(1/v_{\mathbf{M}})'(t)}{-\log(v_{\mathbf{M}}(t))}>1,
\end{equation}
which obviously implies
\begin{equation}\label{Thm39liminf0}
\liminf_{t\rightarrow+\infty}\frac{tv_{\mathbf{M}}(t)(1/v_{\mathbf{M}})'(t)}{-\log(v_{\mathbf{M}}(t))}>0.
\end{equation}

We characterize now the weaker property \eqref{Thm39liminf0} in terms of a known growth condition for $\mathbf{M}$; before proving this technical result we note:

\begin{remark}\label{omegaMvanishingrem}
By definition \eqref{weights} one has $-\log(v_{\mathbf{M}}(t))=\omega_{\mathbf{M}}(t)$ for all $t\ge 0$, but $\omega_{\mathbf{M}}(t)=0$ for all $t\in[0,\mu_1]$ which follows by Lemma \ref{lemma1} resp. by \eqref{assointrepr}. Therefore, whenever considering $\omega_{\mathbf{M}}(t)$ in the denominator of a fraction it is natural to investigate (all) $t>\mu_1$ which is not restricting the generality since for $t\mapsto\omega_{\mathbf{M}}(t)$ one is interested in having $t\rightarrow+\infty$ and thus the relevant information is ``localized at infinity.'' For example, see the proof of \cite[Prop. A.1 $(ii)$]{conjugateweightfunction} where we restricted to all integers $p\ge 2$ sufficiently large resp. to all $t>\mu_1$ and in the first part of \cite[Prop. A.1 $(i)$]{conjugateweightfunction} one should also consider $t>\mu_1$ but this issue does not effect the validity of this statement.
\end{remark}

\begin{lemma}\label{Thm39Lemma}
Let $\mathbf{M}$ be a log-convex weight sequence, then the following are equivalent:
\begin{itemize}
\item[$(i$)] \eqref{Thm39liminf0} is valid.

\item[$(ii)$] $\mathbf{M}$ satisfies \hyperlink{mg}{$(\on{mg})$}.	
\end{itemize}
\end{lemma}

\emph{Note:} \hyperlink{mg}{$(\on{mg})$} is obviously preserved under equivalence of (weight) sequences and this result extends the known characterization for \hyperlink{mg}{$(\on{mg})$} (see \cite{modgrowthstrange} and the citations therein).

\demo{Proof}
$(i)\Rightarrow(ii)$ First, let us compute the expression in \eqref{Thm39liminf0} explicitly: If $\mu_p\le t<\mu_{p+1}$ with $\mu_p>\mu_1$ (see Remark \ref{omegaMvanishingrem}), then similarly as before by definition and Lemma \ref{lemma1} one has
\begin{equation}\label{Thm39Lemmaequ}
\frac{t\exp(-\omega_{\mathbf{M}}(t))\exp(\omega_{\mathbf{M}}(t))\omega_{\mathbf{M}}'(t)}{\omega_{\mathbf{M}}(t)}=\frac{t(p/t)}{\omega_{\mathbf{M}}(t)}=\frac{p}{\omega_{\mathbf{M}}(t)}=\frac{p}{\log(t^p/M_p)}.
\end{equation}
If \eqref{Thm39liminf0} holds, then there exist $\epsilon>0$ and $p_0\in\NN_{>0}$ such that for all $t\ge\mu_{p_0}$ one has by \eqref{Thm39Lemmaequ} that $\frac{p}{\log(t^p/M_p)}>\epsilon\Leftrightarrow t<\exp(\epsilon^{-1})(M_p)^{1/p}$ if $\mu_p\le t<\mu_{p+1}$ and $p\ge p_0$. For this recall that $\lim_{p\rightarrow+\infty}\mu_p=+\infty$ holds and w.l.o.g. one can assume $\mu_{p_0}<\mu_{p_0+1}$. The particular choice $t=\mu_p$ gives $\mu_p<\exp(\epsilon^{-1})(M_p)^{1/p}$ and so \eqref{quotientrootcompprop} is shown for all $p\ge p_0$ such that $\mu_p<\mu_{p+1}$ when taking $A:=\exp(\epsilon^{-1})$. Moreover, the above estimate is strict and hence even gives, by considering the limit $t=\mu_{p+1}$, that $\mu_{p+1}\le\exp(\epsilon^{-1})(M_p)^{1/p}$ for all $p\ge p_0$ such that $\mu_p<\mu_{p+1}$.

If now $p\ge p_0$ and $\mu_p<\mu_{p+1}=\dots=\mu_{p+\ell+1}<\mu_{p+\ell+2}$ for some $\ell\in\NN_{>0}$, then by the above and the fact that $p\mapsto(M_p)^{1/p}$ is non-decreasing (by log-convexity and $M_0=1$) we obtain  $\mu_{p+i}=\mu_{p+1}\le\exp(\epsilon^{-1})(M_{p})^{1/p}\le\exp(\epsilon^{-1})(M_{p+i})^{1/(p+i)}$ for all $1\le i\le\ell+1$. Thus \eqref{quotientrootcompprop} is shown for all $p\ge p_0$ again with $A:=\exp(\epsilon^{-1})$. And by enlarging $A$ sufficiently (if necessary) we get \eqref{quotientrootcompprop} for all $p\in\NN_{>0}$ and hence \hyperlink{mg}{$(\on{mg})$}.\vspace{6pt}

$(ii)\Rightarrow(i)$ As said before, condition \hyperlink{mg}{$(\on{mg})$} is equivalent to $\mu_p\le A(M_p)^{1/p}$ for some $A>1$ (w.l.o.g.) and all $p\in\NN_{>0}$. Then for any $t\ge\mu_1$ we can find $p\in\NN_{>0}$ such that $\mu_p\le t<\mu_{p+1}$ and one gets $\omega_{\mathbf{M}}(t)\le\omega_{\mathbf{M}}(\mu_{p+1})=\log\left(\frac{\mu_{p+1}^{p+1}}{M_{p+1}}\right)\le(p+1)\log(A)\le 2p\log(A)$. This estimate gives $\frac{p}{\omega_{\mathbf{M}}(t)}\ge\frac{1}{2\log(A)}>0$ for any $t>\mu_1$. (However, by the convention $\frac{p}{0}=+\infty$, this estimate is formally valid even for all $t\ge 0$.) Thus, in view of \eqref{Thm39Lemmaequ}, property \eqref{Thm39liminf0} is shown.
\qed\enddemo

Therefore, Lemma \ref{Thm39Lemma} implies the fact that for any log-convex weight sequence violating \hyperlink{mg}{$(\on{mg})$} one cannot apply \cite[Thm. 39]{Bonet2022survey} dealing with the description of the proper closed invariant subspaces for the (by assumption) continuous integration operator $J: H^{0}_{v_{\mathbf{M}}}(\CC)\rightarrow H^{0}_{v_{\mathbf{M}}}(\CC)$: the crucial (but only sufficient) $\liminf$-condition \eqref{Thm39liminf} in this result fails.

\subsection{The unit disk case}\label{unitdiscsection}
Recall \eqref{weightfctunitdisk}, see also \cite[Sect. 7]{solidassociatedweight} (with $c=1$), and hence we consider the weight
$$v_{\mathbf{M},\mathbb{D}}(t):=\exp\left(-\omega_{\mathbf{M}}\left(\frac{1}{1-t}\right)\right),\;\;\;t\in[0,1).$$
$v_{\mathbf{M},\mathbb{D}}$ is continuous, non-increasing and rapidly decreasing which means that $\lim_{t\rightarrow 1^{-}}t^pv_{\mathbf{M},\mathbb{D}}(t)=0$ for all $p\in\NN$ (recall \cite[Sect. 2, p. 3]{Bonet2022survey}). For the following, crucially note (again) the different notation in \cite{AbaninThien17} compared with \cite{Bonet2022survey} and so having $v\mapsto v^{-1}$ when comparing the conditions for the weights in these articles.\vspace{6pt}

Related to \cite[Thm. 27 \& 30, Cor. 31]{Bonet2022survey} resp. \cite[Lemma 2.6 \& Prop. 3.12]{AbaninThien17} is the study of the following expression:
\begin{equation}\label{crucialdiskcondition}
-\frac{(1-t)v'_{\mathbf{M},\mathbb{D}}(t)}{v_{\mathbf{M},\mathbb{D}}(t)}=-\frac{(1-t)(\exp(-\omega_{\mathbf{M},\mathbb{D}}(t)))'}{\exp(-\omega_{\mathbf{M},\mathbb{D}}(t))}.
\end{equation}

Let now $\mathbf{M}$ be a log-convex weight sequence and let $t\in[0,1)$ be such that $1-\frac{1}{\mu_p}\le t<1-\frac{1}{\mu_{p+1}}$, equivalently $\mu_p\le\frac{1}{1-t}<\mu_{p+1}$, for some $p\in\NN_{>0}$. Involving Lemma \ref{lemma1} we get that
\begin{align*}
v_{\mathbf{M},\mathbb{D}}(t)&=\exp\left(-\omega_{\mathbf{M}}\left(\frac{1}{1-t}\right)\right)=\exp\left(-\log\left(\frac{(1/(1-t))^p}{M_p}\right)\right)=\exp\left(\log\left(\frac{M_p}{(1/(1-t))^p}\right)\right)
\\&
=M_p(1-t)^p.
\end{align*}
Consequently, for all such values $t$ we obtain
\begin{equation*}\label{Thm30Cor31cond}
-\frac{(1-t)v'_{\mathbf{M},\mathbb{D}}(t)}{v_{\mathbf{M},\mathbb{D}}(t)}=-\frac{(1-t)M_pp(1-t)^{p-1}(-1)}{M_p(1-t)^p}=p,
\end{equation*}
and as $t\rightarrow 1^{-}$ we get that $p\rightarrow+\infty$ and hence
\begin{equation}\label{Thm30Cor31condcons}
\lim_{t\rightarrow 1^{-}}\left(-\frac{(1-t)v'_{\mathbf{M},\mathbb{D}}(t)}{v_{\mathbf{M},\mathbb{D}}(t)}\right)=+\infty.
\end{equation}
Thus the next result is shown:

\begin{theorem}\label{Thm30&31}
Let $\mathbf{M}$ be a log-convex weight sequence and put $w_{\mathbf{M},\mathbb{D}}(t):=(1-t)v_{\mathbf{M},\mathbb{D}}(t)$, $t\in[0,1)$. Then the following hold:
\begin{itemize}
\item[$(i)$] The integration operator $J: H^{\infty}_{w_{\mathbf{M},\mathbb{D}}}(\mathbb{D})\rightarrow H^{\infty}_{v_{\mathbf{M},\mathbb{D}}}(\mathbb{D})$ is continuous.

\item[$(ii)$] The differentiation operator $D: H^{\infty}_{v_{\mathbf{M},\mathbb{D}}}(\mathbb{D})\rightarrow H^{\infty}_{w_{\mathbf{M},\mathbb{D}}}(\mathbb{D})$ is not continuous and surjective (i.e. $D$ is either not continuous or not surjective).

\item[$(iii)$] The differentiation operator $D: H^{\infty}_{v_{\mathbf{M},\mathbb{D}}}(\mathbb{D})\rightarrow H^{\infty}_{w_{\mathbf{M},\mathbb{D}}}(\mathbb{D})$ is not continuous.
\end{itemize}
\end{theorem}

\emph{Note:} In \cite[Thm. 30 \& Cor. 31]{Bonet2022survey} one should replace the weight $w(t):=\frac{v(t)}{1-t}$ by $w(t):=v(t)(1-t)$ and hence use the same expression as in \cite[Thm. 27]{Bonet2022survey}. Moreover, by \cite[Lemma 3.15, Cor. 3.16, Rem. 3.17]{AbaninThien17} a full characterization of the continuity of $J$ in terms of the growth of the expression in \eqref{crucialdiskcondition} is unclear and only the sufficiency is shown.

\demo{Proof}
First, recall that by Lemma \ref{disclogconvexlemma} above the weight $v_{\mathbf{M},\mathbb{D}}$ is radial and (log-)convex. Then $(i)$ follows from \eqref{Thm30Cor31condcons} and \cite[Thm. 30 (iii)$\Rightarrow$(i)]{Bonet2022survey} resp. \cite[Lemma 3.15 \& Cor. 3.16]{AbaninThien17}. For $(ii)$ we use \eqref{Thm30Cor31condcons} and \cite[Cor. 31]{Bonet2022survey} resp. \cite[Cor. 3.18]{AbaninThien17}; indeed \eqref{Thm30Cor31condcons} violates the second half of the characterizing condition in these results. Finally, $(iii)$ is valid by \eqref{Thm30Cor31condcons} and \cite[Thm. 27 (i)$\Leftrightarrow$(iii)]{Bonet2022survey} resp. \cite[Lemma 2.6 \& Thm. 2.8]{AbaninThien17}.
\qed\enddemo

For reasons of consistency we investigate now some more technical properties on weights appearing in the characterizing statements from \cite{Bonet2022survey} and \cite{AbaninThien17} mentioned in the proof before. The forthcoming also illustrates the fact that the weight sequence setting is rigid in the unit disk case in the sense that in the alternative characterizing conditions the limits exist and attain extreme values, i.e. coincide with $0$ resp. $+\infty$. In particular, no oscillation of the crucial quotients is possible.

In \cite[Thm. 27 (iv)]{Bonet2022survey} resp. \cite[Lemma 2.6 $(iv)$]{AbaninThien17}, for a weight function $v$ on $\mathbb{D}$ the next condition appears
\begin{equation}\label{Thm274}
\sup_{n\in\NN}\frac{v(1-2^{-n})}{v(1-2^{-n-1})}<+\infty,
\end{equation}
whereas in \cite[Thm. 30 (iv)]{Bonet2022survey} resp. \cite[Lemma 3.15 $(iv)$]{AbaninThien17} one considers
\begin{equation}\label{Thm304}
\exists\;k\in\NN_{>0}:\;\;\;\limsup_{n\rightarrow+\infty}\frac{v(1-2^{-n-k})}{v(1-2^{-n})}<1.
\end{equation}
Note that in \cite[Thm. 27 (iv)]{Bonet2022survey} there is a misprint and one should read $v(1-2^{-n-1})$ instead of $v(1-2^{-n+1})$ in the denominator. Moreover, in \cite[Thm. 27 (v)]{Bonet2022survey} resp. \cite[Lemma 2.6 $(vi)$]{AbaninThien17} condition
\begin{equation}\label{thm27v}
v(t^2)=O(v(t)),\;\;\;t\rightarrow 1^{-},
\end{equation}
appears. Finally, due to the different notation $v$ vs. $v^{-1}$ and the comparison with \cite[Lemma 3.15 $(vi)$]{AbaninThien17}, the condition in \cite[Thm. 30 (v)]{Bonet2022survey} should be replaced by
\begin{equation}\label{thm30v}
\exists\;\gamma>1:\;\;\;\limsup_{t\rightarrow 1^{-}}\frac{v(t)}{v(t^{\gamma})}<1.
\end{equation}

\begin{lemma}\label{Thm30&31lemma}
Let $\mathbf{M}$ be a log-convex weight sequence, then
\begin{equation}\label{Thm30&31lemmaequ}
\forall\;k\in\NN_{>0}:\;\;\;\lim_{n\rightarrow+\infty}\frac{v_{\mathbf{M},\mathbb{D}}(1-2^{-n-k})}{v_{\mathbf{M},\mathbb{D}}(1-2^{-n})}=0,
\end{equation}
and
\begin{equation}\label{Thm30&31lemmaequ2}
\forall\;\gamma>1:\;\;\;\lim_{t\rightarrow 1^{-}}\frac{v_{\mathbf{M},\mathbb{D}}(t^{\gamma})}{v_{\mathbf{M},\mathbb{D}}(t)}=+\infty.
\end{equation}
\end{lemma}

\emph{We summarize:}

\begin{itemize}
\item[$(*)$] \eqref{Thm30&31lemmaequ2} yields the fact that \eqref{thm27v} fails but \eqref{thm30v} holds even for any $\gamma>1$. Moreover, \eqref{Thm30&31lemmaequ} contradicts \eqref{Thm274} and implies \eqref{Thm304}.

\item[$(*)$] This is consistent with the proof of $(i)$ resp. $(iii)$ in Theorem \ref{Thm30&31}, where we involved \eqref{Thm30Cor31condcons}, and with \cite[Thm. 30 (iii)$\Leftrightarrow$(iv)$\Leftrightarrow$(v)]{Bonet2022survey} (the condition in (v) being replaced by \eqref{thm30v}) resp. with \cite[Thm. 27 (iii)$\Leftrightarrow$(iv)$\Leftrightarrow$(v)]{Bonet2022survey} and cf. also \cite[Lemma 3.15 $(i)\Leftrightarrow(iv)\Leftrightarrow(vi)$]{AbaninThien17} resp. \cite[Lemma 2.6 $(i)\Leftrightarrow(iv)\Leftrightarrow(vi)$]{AbaninThien17}.
\end{itemize}

\demo{Proof}
By definition \eqref{weightfctunitdisk} one has $$v_{\mathbf{M},\mathbb{D}}(1-2^{-n-k})=\exp\left(-\omega_{\mathbf{M}}\left(\frac{1}{1-(1-2^{-n-k})}\right)\right)=\exp(-\omega_{\mathbf{M}}(2^{n+k})),$$ and similarly $v_{\mathbf{M},\mathbb{D}}(1-2^{-n})=\exp(-\omega_{\mathbf{M}}(2^{n}))$. Thus, for any fixed $k\in\NN_{>0}$ it holds that
$$\frac{v_{\mathbf{M},\mathbb{D}}(1-2^{-n-k})}{v_{\mathbf{M},\mathbb{D}}(1-2^{-n})}=\exp(\omega_{\mathbf{M}}(2^{n})-\omega_{\mathbf{M}}(2^{n+k}))\rightarrow 0,\;\;\;n\rightarrow+\infty,$$
since $\omega_{\mathbf{M}}(at)-\omega_{\mathbf{M}}(t)\rightarrow+\infty$ as $t\rightarrow+\infty$ for any (fixed) $a>1$ (see \cite[Lemma 2.9]{weightedentireinclusion1}) and take $a:=2^k$. Thus \eqref{Thm30&31lemmaequ} is shown.\vspace{6pt}

Now we prove \eqref{Thm30&31lemmaequ2} and so let $\gamma>1$ be arbitrary but from now on fixed. First, similarly as before by \eqref{weightfctunitdisk} we get $$\frac{v_{\mathbf{M},\mathbb{D}}(t^{\gamma})}{v_{\mathbf{M},\mathbb{D}}(t)}=\exp\left(\omega_{\mathbf{M}}\left(\frac{1}{1-t}\right)-\omega_{\mathbf{M}}\left(\frac{1}{1-t^{\gamma}}\right)\right)=\exp\left(\omega_{\mathbf{M}}(s)-\omega_{\mathbf{M}}\left(\frac{s^{\gamma}}{s^{\gamma}-(s-1)^{\gamma}}\right)\right),$$
with the substitution $s:=\frac{1}{1-t}\Leftrightarrow t=1-\frac{1}{s}$, $t\in[0,1)$, and so $s\in[1,+\infty)$. For this note that $t^{\gamma}=\frac{(s-1)^{\gamma}}{s^{\gamma}}$ and $\frac{1}{1-t^{\gamma}}=\frac{s^{\gamma}}{s^{\gamma}-(s-1)^{\gamma}}$ is valid.
Now select $a\in(1,\gamma)$ and note that $s^{\gamma}-(s-1)^{\gamma}=\gamma\xi^{\gamma-1}$ for some $\xi\in(s-1,s)$ which implies that $\gamma s^{\gamma-1}\ge s^{\gamma}-(s-1)^{\gamma}\ge\gamma(s-1)^{\gamma-1}$ for all $s\in[1,+\infty)$. In particular, this also proves that $\lim_{s\rightarrow+\infty}\frac{s^{\gamma}}{s^{\gamma}-(s-1)^{\gamma}}=+\infty$ and note that $r\rightarrow 1^{-}$ if and only if $s\rightarrow+\infty$. Thus, if $s>1$, then $a\frac{s^{\gamma}}{s^{\gamma}-(s-1)^{\gamma}}\le\frac{a}{\gamma}\frac{s^{\gamma}}{(s-1)^{\gamma-1}}$ and we estimate as follows:
\begin{equation}\label{Thm30&31lemmaaux}
\frac{a}{\gamma}\frac{s^{\gamma}}{(s-1)^{\gamma-1}}\le s\Leftrightarrow\frac{s}{s-1}\le\left(\frac{\gamma}{a}\right)^{1/(\gamma-1)}.
\end{equation}
Since $\frac{\gamma}{a}>1$, $\gamma-1>0$ and since the left-hand side in \eqref{Thm30&31lemmaaux} tends to $1$ as $s\rightarrow+\infty$ this estimate is valid for all $s$ sufficiently large and so we have shown that
$$\exists\;s_{\gamma}>1\;\forall\;s\ge s_{\gamma}:\;\;\;a\frac{s^{\gamma}}{s^{\gamma}-(s-1)^{\gamma}}\le s.$$
Therefore, for all $s\ge s_{\gamma}$ we infer
\begin{equation}\label{Thm30&31lemmaequ3}
\omega_{\mathbf{M}}(s)-\omega_{\mathbf{M}}\left(\frac{s^{\gamma}}{s^{\gamma}-(s-1)^{\gamma}}\right)\ge\omega_{\mathbf{M}}\left(a\frac{s^{\gamma}}{s^{\gamma}-(s-1)^{\gamma}}\right)-\omega_{\mathbf{M}}\left(\frac{s^{\gamma}}{s^{\gamma}-(s-1)^{\gamma}}\right),
\end{equation}
and by \cite[Lemma 2.9]{weightedentireinclusion1} applied to chosen $a\in(1,\gamma)$ the right-hand side in \eqref{Thm30&31lemmaequ3} tends to $+\infty$ as $s\rightarrow+\infty$. Thus \eqref{Thm30&31lemmaequ2} follows.
\qed\enddemo

\bibliographystyle{plain}
\bibliography{Bibliography}
\end{document}